\documentstyle[12pt]{article}
\newlength{\dinwidth}
\newlength{\dinmargin}
\newtheorem{Definition}{Definition}[section]
\newtheorem{Theorem}[Definition]{Theorem}
\newtheorem{Proposition}[Definition]{Proposition}
\newtheorem{Lemma}[Definition]{Lemma}

\def\AA{{\cal A}}
\def\lcrc{{\mbox{\footnotesize $\circ$}} }
\def\CC{{\bf C}}
\def\RR{{\bf R}}
\def\NN{{\bf N}}
\def\qS{{\sf q}(S,\sigma)}
\def\prS{{\sf pr}(S,\sigma)}
\def\puS{{\sf pu}(S,\sigma)}
\def\Um{U_{\mu}}
\def\RmB{|R_{\mu}|}
\def\AS{{\cal A}[S,\sigma]}
\def\om{\omega_{\mu}}
\def\Hh{{\cal H}}
\def\Ff{{\cal F}}
\def\mT{\tilde{\mu}}
\def\MT{\tilde{M}}
\def\ST{\tilde{\Sigma}}
\def\gT{\tilde{\gamma}}
\def\Rr{{\cal R}}
\def\AO{{\cal A(O)}}
\def\Ss{{\cal S}}
\def\Oo{{\cal O}}
\def\DS{{\cal D}_{\Sigma}}
\def\CoS{C_0^{\infty}(\Sigma,{\bf R})}
\def\CoM{C_0^{\infty}(M,{\bf R})}
\def\uu{u_0 \oplus u_1}
\def\vv{v_0 \oplus v_1}
\begin{document}
\title{Continuity of Symplectically Adjoint Maps\\
and the Algebraic Structure \\  of Hadamard Vacuum Representations
\\ for Quantum Fields on Curved Spacetime}
 \author{{\sc Rainer Verch}\thanks{Supported by
a Von Neumann Fellowship of the Operator Algebras Network,
EC Human Capital and Mobility Programme.}
\\[12pt]
\normalsize
        Dipartimento di Matematica \\
\normalsize        Universit\`a di Roma II ``Tor Vergata''\\
\normalsize        I - 00133 Roma, Italy \\
\normalsize    e-mail: verch$@$x4u2.desy.de}
  \date{\normalsize Revised version/ 12 November 1996 }
\maketitle
${}$
\\[24pt] {\small
{\bf Abstract.} We derive for a pair of operators on a symplectic
space which are adjoints of each other with respect to the symplectic
form (that is, they are sympletically adjoint) that, if they are
bounded for some scalar product on the symplectic space dominating
the symplectic form, then they are bounded with respect to a
one-parametric family of scalar products canonically associated
with the initially given one, among them being its ``purification''.
As a typical example we consider a scalar field on a globally
hyperbolic spacetime governed by the Klein-Gordon equation; the classical
system is described by a symplectic space and the temporal evolution
by symplectomorphisms (which are symplectically adjoint to
their inverses). A
natural scalar product is that inducing the classical energy norm,
and an application of the above result yields that its ``purification''
induces on the one-particle space of the quantized system a
topology which coincides with that given by
 the two-point functions of quasifree Hadamard states.
These findings will be shown to lead to new results concerning the
structure of the local (von Neumann) observable-algebras in representations
of quasifree Hadamard states of the Klein-Gordon field in an arbitrary
globally hyperbolic spacetime, such as local definiteness, local
primarity and Haag-duality (and also split- and type ${\rm III}_1$-properties).
A brief review of this circle of notions, as well as of properties of
Hadamard states, forms part of the article. } 
\noindent
\section{Introduction}
In the first part of this paper we shall investigate a
special case of relative continuity of symplectically adjoint maps
of a symplectic space. By this, we mean the following.
Suppose that $(S,\sigma)$ is a symplectic space, i.e.\  $S$
is a real-linear vector space with an anti-symmetric,
non-degenerate bilinear form $\sigma$ (the symplectic form).
A pair $V,W$ of linear maps of $S$ will be called
{\it symplectically adjoint} if
$\sigma(V\phi,\psi) = \sigma(\phi,W\psi)$ for all $\phi,\psi \in S$.
 Let $\mu$ and $\mu'$ be two
scalar products on $S$  and assume that,
for each pair $V,W$ of symplectically adjoint
linear maps of $(S,\sigma)$, the boundedness
of both $V$ and $W$ with respect to $\mu$
 implies their boundedness with respect to $\mu'$.
Such a situation we refer to as {\it relative $\mu - \mu'$
continuity of symplectically adjoint maps} (of $(S,\sigma)$).
A particular example of symplectically adjoint maps is
provided by the pair $T,T^{-1}$ whenever $T$ is a symplectomorphism
of $(S,\sigma)$. (Recall that a symplectomorphism of $(S,\sigma)$
is a bijective linear map $T : S \to S$ which preserves the
symplectic form, $\sigma(T\phi,T\psi) = \sigma(\phi,\psi)$ for
all $\phi,\psi \in S$.)

In the more specialized case to be considered in the present work,
which will soon be indicated to be relevant in applications,
we show that a certain distinguished relation between a
 scalar product $\mu$ on $S$ 
 and a second one, $\mu'$, 
 is sufficient for the relative $\mu - \mu'$
continuity of symplectically adjoint maps.
(We give further details in Chapter 2,
and in the next paragraph.) 
 The result will be applied in Chapter 3 to answer
a couple of open questions
concerning the algebraic structure of the quantum theory of the free
scalar field in arbitrary globally hyperbolic spacetimes:
the local definiteness, local primarity  and Haag-duality
in representations of the local observable algebras
induced by quasifree Hadamard states, 
as well as the determination of the type of the local
von Neumann algebras in such representations.
Technically, what needs to be proved in our approach to this problem
is the continuity of the temporal evolution of the Cauchy-data of
solutions of the scalar Klein-Gordon equation
\begin{equation}
(\nabla^a \nabla_a + r)\varphi = 0
\end{equation}
in a globally hyperbolic spacetime with respect to a certain
topology on the Cauchy-data space.
(Here, $\nabla$ is the covariant derivative of the metric $g$
on the spacetime, and $r$ an arbitrary realvalued,
 smooth function.)
The Cauchy-data space is a symplectic space on which the said
temporal evolution is realized by symplectomorphisms. It
turns out that the classical ``energy-norm'' of solutions
of (1.1), which is given by  a scalar
product $\mu_0$ on the Cauchy-data space, and the
topology relevant for the required continuity statement
(the ``Hadamard one-particle space norm''), induced by a
scalar product $\mu_1$ on the Cauchy-data space, are precisely
in the relation for which our result on relative $\mu_0 - \mu_1$
continuity of symplectically adjoint maps applies. Since the continuity
of the Cauchy-data evolution in the classical energy norm,
i.e.\ $\mu_0$, is well-known, the desired continuity in the
$\mu_1$-topology follows.

The argument just described may be viewed as the prime example of
application of the relative continuity result. In fact,
the relation between $\mu_0$ and $\mu_1$ is abstracted from
the relation between the classical energy-norm and the
one-particle space norms arising from ``frequency-splitting'' procedures
in the canonical quantization of (linear) fields.
 This relation has been made precise
in a recent paper by Chmielowski [11]. It provides the
starting point for our investigation in Chapter 2, where
we shall see  that one can associate
with a dominating scalar product $\mu \equiv \mu_0$ on
$S$ in a canonical way a positive, symmetric operator
$\RmB$ on the $\mu$-completion of $S$, and a family of scalar
products $\mu_s$, $s > 0$, on $S$, defined as $\mu$ with
$\RmB^s$ as an operator kernel. Using abstract
interpolation, it will be shown that
then relative $\mu_0 - \mu_s$ continuity of symplectically adjoint maps
holds for all $0 \leq s \leq 2$. The relative
$\mu_0 - \mu_1$ continuity arises as
a special case.
In fact, it turns out that the indicated interpolation
argument may even be extended to an apparently more general
situation from which the relative $\mu_0 - \mu_s$ continuity
of symplectically adjoint maps derives as a corollary, see
Theorem 2.2.

Chapter 3 will be concerned with the application of the result
of Thm.\ 2.2 as indicated above. In the preparatory Section
3.1, some notions of general relativity will be summarized, along
with the introduction of some notation. Section 3.2 contains a brief
synopsis of the notions of local definiteness, local primarity and
Haag-duality in the the context of quantum field theory in curved
spacetime. In Section 3.3 we present the $C^*$-algebraic quantization of the
KG-field obeying (1.1) on a globally hyperbolic spacetime, following [16].
Quasifree Hadamard states will be described in Section 3.4 according
to the definition given in [45]. In the same section we briefly summarize
some properties of Hadamard two-point functions, and derive, in
Proposition 3.5, the result concerning the continuity of the
Cauchy-data evolution maps in the topology of the Hadamard two-point
functions which was mentioned above. It will be seen in the last Section
3.5 that this leads, in combination with results obtained earlier
[64,65,66], to Theorem 3.6 establishing detailed properties of the algebraic
structure of the local von Neumann observable algebras in representations
induced by quasifree Hadamard states of the Klein-Gordon field over
an arbitrary globally hyperbolic spacetime.
\section{Relative Continuity of Symplectically Adjoint Maps}
\setcounter{equation}{0}
Let $(S,\sigma)$ be a symplectic space. A (real-linear) scalar
product $\mu$ on $S$ is said to {\it dominate} $\sigma$
if the estimate
\begin{equation}
|\sigma(\phi,\psi)|^2 \leq 4 \cdot \mu(\phi,\phi)\,\mu(\psi,\psi)\,,
\quad \phi,\psi \in S\,,
\end{equation}
holds; the set of all scalar products on $S$ which dominate $\sigma$
will be denoted by ${\sf q}(S,\sigma)$.
Given $\mu \in {\sf q}(S,\sigma)$, we write $H_{\mu} \equiv
\overline{S}^{\mu}$ for the completion of $S$ with respect to the
topology induced by $\mu$, and denote by $\sigma_{\mu}$ the
$\mu$-continuous extension, guaranteed to uniquely exist by (2.1),
of $\sigma$ to $H_{\mu}$. The estimate (2.1) then extends to
$\sigma_{\mu}$ and all $\phi,\psi \in H_{\mu}$. This entails
that there is a uniquely determined, $\mu$-bounded linear
operator $R_{\mu} : H_{\mu} \to H_{\mu}$ with the property
\begin{equation}
\sigma_{\mu}(x,y) = 2\,\mu(x,R_{\mu}y)\,, \quad x,y \in H_{\mu}\,.
\end{equation}
The antisymmetry of $\sigma_{\mu}$ entails for the
$\mu$-adjoint $R_{\mu}^*$ of $R_{\mu}$
\begin{equation}
R_{\mu}^* = - R_{\mu}\,,
\end{equation}
and by (2.1) one finds that the operator norm of $R_{\mu}$
is bounded by 1, $||\,R_{\mu}\,|| \leq 1$.
The operator $R_{\mu}$ will be called the {\it polarizator} of $\mu$.

In passing, two things should be noticed here:
\\[6pt]
(1) $R_{\mu}|S$ is injective since $\sigma$ is a non-degenerate
bilinear form on $S$, but $R_{\mu}$ need not be injective on
on all of $H_{\mu}$, as $\sigma_{\mu}$ may be degenerate.
\\[6pt]
(2) In general, it is not the case that $R_{\mu}(S) \subset S$.
\\[6pt]
Further properties of $R_{\mu}$ will be explored below.
Let us first focus on two significant subsets of $\qS$ which
are intrinsically characterized by properties of the
corresponding $\sigma_{\mu}$ or, equivalently, the $R_{\mu}$.

The first is $\prS$, called the set of {\it primary}
scalar products on $(S,\sigma)$, where $\mu \in \qS$ is
in $\prS$ if $\sigma_{\mu}$ is a symplectic form
(i.e.\ non-degenerate) on $H_{\mu}$. In view of (2.2) and
(2.3), one can see that this is equivalent to either
(and hence, both) of the following conditions:
\begin{itemize}
\item[(i)] \quad $R_{\mu}$ is injective,
\item[(ii)] \quad $R_{\mu}(H_{\mu})$ is dense in $H_{\mu}$.
\end{itemize}
The second important subset of $\qS$ is denoted by
$\puS$ and defined as consisting of those $\mu \in \qS$
which satisfy the {\it saturation property}
\begin{equation}
\mu(\phi,\phi) = \sup_{\psi \in S\backslash \{0\} } \,
 \frac{|\sigma(\phi,\psi)|^2}{4 \mu(\psi,\psi) } \,,\ \ \  \psi \in S \,.
\end{equation}
The set $\puS$ will be called the set of {\it pure} scalar
products on $(S,\sigma)$. It is straightforward to check that
$\mu \in \puS$ if and only if $R_{\mu}$ is a unitary
anti-involution, or complex structure, i.e.\
$R_{\mu}^{-1} = R_{\mu}^*$, $R_{\mu}^2 = - 1$. Hence
$\puS \subset \prS$.
\\[10pt]
Our terminology reflects well-known relations between properties of
quasifree states on the (CCR-) Weyl-algebra of a symplectic space
$(S,\sigma)$ and properties of $\sigma$-dominating scalar products
on $S$, which we shall briefly recapitulate. We refer to
[1,3,5,45,49] and also references quoted therein for proofs
and further discussion of the following statements.

Given a symplectic space $(S,\sigma)$, one can associate with it
uniquely (up to $C^*$-algebraic equivalence) a $C^*$-algebra
$\AS$, which is generated by a family of unitary elements
$W(\phi)$, $\phi \in S$, satisfying the canonical commutation
relations (CCR) in exponentiated form,
\begin{equation}
W(\phi)W(\psi) = {\rm e}^{-i\sigma(\phi,\psi)/2}W(\phi + \psi)\,,
\quad \phi,\psi \in S\,.
\end{equation}
The algebra $\AS$ is called the {\it Weyl-algebra}, or
{\it CCR-algebra}, of $(S,\sigma)$. It is not difficult to see that
if $\mu \in \qS$, then one can define a state (i.e., a positive,
normalized linear functional) $\om$ on $\AS$ by setting
\begin{equation}
 \om (W(\phi)) : = {\rm e}^{- \mu(\phi,\phi)/2}\,, \quad \phi \in S\,.
\end{equation}
Any state on the Weyl-algebra $\AS$ which can be realized in this way
is called a {\it quasifree state}. Conversely, given any quasifree
state $\om$ on $\AS$, one can recover its $\mu \in \qS$ as
\begin{equation}
\mu(\phi,\psi) = 2 {\sf Re}\left. \frac{\partial}{\partial t}
\frac{\partial}{\partial \tau} \right|_{t = \tau = 0}
 \om (W(t\phi)W(\tau \psi))\,, \quad \phi,\psi \in S\,.
\end{equation}
So there is a one-to-one correspondence between quasifree states
on $\AS$ and dominating scalar products on $(S,\sigma)$.
\\[10pt]
Let us now recall the subsequent terminology. To a state $\omega$
on a $C^*$-algebra $\cal B$ there corresponds (uniquely up to
unitary equivalence) a triple $(\Hh_{\omega},\pi_{\omega},\Omega_{\omega})$,
called the GNS-representation of $\omega$ (see e.g.\ [5]), characterized by
the following properties: $\Hh_{\omega}$ is a complex Hilbertspace,
$\pi_{\omega}$ is a representation of $\cal B$ by bounded linear operators
on $\Hh_{\omega}$ with cyclic vector $\Omega_{\omega}$, and
$\omega(B) = \langle \Omega_{\omega},\pi_{\omega}(B)\Omega_{\omega}
\rangle $ for all  $B \in \cal B$.
Hence one is led to associate with $\omega$ and
$\cal B$ naturally
        the $\omega$-induced von Neumann algebra $\pi_{\omega}({\cal B})^-$,
where the bar means taking the closure with respect to the weak
operator topology in the set of bounded linear operators on $\Hh_{\omega}$.
One refers to $\omega$ (resp., $\pi_{\omega}$) as {\it primary}
if $\pi_{\omega}({\cal B})^- \cap \pi_{\omega}({\cal B})' = \CC \cdot 1$
(so the center of $\pi_{\omega}({\cal B})^-$ is trivial), where the
prime denotes taking the commutant, and as {\it pure} if
$\pi_{\omega}({\cal B})' = \CC\cdot 1$ (i.e.\ $\pi_{\omega}$ is
irreducible --- this is equivalent to the statement that $\omega$
is not a (non-trivial) convex sum of different states).

In the case where $\om$ is a quasifree state on a Weyl-algebra
$\AS$, it is known that (cf.\ [1,49])
\begin{itemize}
\item[(I)] $\om$ is primary if and only if $\mu \in \prS$,
\item[(II)] $\om$ is pure if and only if $\mu \in \puS$.
\end{itemize}
${}$\\
We return to the investigation of the properties of the polarizator
$R_{\mu}$ for a dominating scalar product $\mu$ on a symplectic space
$(S,\sigma)$. It possesses a polar decomposition
\begin{equation}
R_{\mu} = \Um \RmB
\end{equation}
on the Hilbertspace $(H_{\mu},\mu)$, where $\Um$ is an isometry
and $\RmB$ is symmetric and has non-negative spectrum. Since
$R_{\mu}^* = - R_{\mu}$, $R_{\mu}$ is normal and thus
$\RmB$ and $\Um$ commute. Moreover, one has
$\RmB \Um^* = - \Um \RmB$, and hence $\RmB$ and $\Um^*$ commute
as well. One readily observes that $(\Um^* + \Um)\RmB = 0$.
The commutativity can by the spectral calculus be generalized to the
statement that, whenever $f$ is a real-valued, continuous function
on the real line, then
\begin{equation}
 [f(\RmB),\Um] = 0 = [f(\RmB),\Um^*] \,,
\end{equation}
where the brackets denote the commutator.

In a recent work [11], Chmielowski noticed that if one defines
for $\mu \in \qS$ the bilinear form
\begin{equation}
 \mT(\phi,\psi) := \mu (\phi,\RmB \psi)\,, \quad \phi,\psi \in S,
\end{equation}
then it holds that $\mT \in \puS$. The proof of this is straightforward.
That $\mT$ dominates $\sigma$ will be seen in Proposition 2.1 below.
To check the saturation property (2.4) for $\mT$, it suffices to
observe that for given $\phi \in H_{\mu}$, the inequality in the
following chain of expressions:
\begin{eqnarray*}
 \frac{1}{4} | \sigma_{\mu}(\phi,\psi) |^2 & = & |\mu(\phi,\Um \RmB \psi)|^2
 \ = \ |\mu(\phi,-\Um^*\RmB\psi) |^2  \\
 & = & |\mu(\RmB^{1/2}\Um\phi,\RmB^{1/2}\psi)|^2
 \\
 & \leq &  \mu(\RmB^{1/2}\Um\phi,\RmB^{1/2}\Um\phi) \cdot \mu(\RmB^{1/2}\psi,
           \RmB^{1/2}\psi) \nonumber
\end{eqnarray*}
is saturated and becomes an equality upon choosing $\psi \in H_{\mu}$
so that $\RmB^{1/2}\psi$ is parallel to $\RmB^{1/2}\Um \phi$.
Therefore one obtains for all $\phi \in S$
\begin{eqnarray*}
 \sup_{\psi \in S\backslash\{0\}}\, \frac{|\sigma(\phi,\psi)|^2}
{4 \mu(\psi,\RmB \psi) } & = & \mu(\RmB^{1/2}\Um \phi,\RmB^{1/2}\Um \phi)
\\
    & = & \mu(\Um\RmB^{1/2}\phi,\Um \RmB^{1/2} \phi)  \\
    & = & \mT(\phi,\phi)\,,
\end{eqnarray*}
which is the required saturation property.

Following Chmielowski, the scalar product $\mT$ on $S$ associated with
$\mu \in \qS$ will be called the {\it purification} of $\mu$.

It appears natural to associate with $\mu \in \qS$ the family $\mu_s$,
$s > 0$, of symmetric bilinear forms on $S$ given by
\begin{equation}
\mu_s(\phi,\psi) := \mu(\phi,\RmB^s \psi)\,, \quad \phi,\psi \in S\,.
\end{equation}
We will use the convention that $\mu_0 = \mu$.
Observe that $\mT = \mu_1$. The subsequent proposition ensues.
\begin{Proposition}
${}$\\[6pt]
(a) $\mu_s$ is a scalar product on $S$ for each $s \geq 0$. \\[6pt]
(b) $\mu_s$ dominates $\sigma$ for $0 \leq s \leq 1$. \\[6pt]
(c) Suppose that there is some $s \in (0,1)$ such that $\mu_s \in \puS$.
Then $\mu_r = \mu_1$ for all $r > 0$. If it is in addition assumed
that $\mu \in \prS$, then it follows that $\mu_r = \mu_1$ for all
$r \geq 0$, i.e.\ in particular $\mu = \mT$. \\[6pt]
(d) If $\mu_s \in \qS$ for some $s > 1$, then $\mu_r = \mu_1$ for
all $r > 0$. Assuming additionally $\mu \in \prS$,  one obtains
$\mu_r = \mu_1$ for all $r \geq 0$, entailing $\mu = \mT$.\\[6pt]
(e) The purifications of the $\mu_s$, $0 < s < 1$, are equal
to $\mT$: We have $\widetilde{\mu_s} = \mT = \mu_1$ for all
$0 < s < 1$.
\end{Proposition}
{\it Proof.} (a) According to (b), $\mu_s$ dominates $\sigma$ for
each $0 \leq s \leq 1$, thus it is a scalar product whenever $s$ is
in that range. However, it is known that
$\mu(\phi,\RmB^s \phi) \geq \mu(\phi,\RmB \phi)^s$ for all vectors
$\phi \in H_{\mu}$ of unit length ($\mu(\phi,\phi) = 1$) and
$1 \leq s < \infty$, cf.\ [60 (p.\ 20)]. This shows that
$\mu_s(\phi,\phi) \neq 0$ for all nonzero $\phi$ in $S$, $s \geq 0$.
\\[6pt]
(b) For $s$ in the indicated range there holds the following estimate:
\begin{eqnarray*}
\frac{1}{4} |\sigma(\phi,\psi)|^2 & = & |\mu(\phi,\Um\RmB \psi)|^2
  \  = \  |\mu(\phi,-\Um^*\RmB \psi )|^2 \\
  & = & | \mu(\RmB^{s/2}U_{\mu} \phi, \RmB^{1 - s/2} \psi) |^2 \\
  & \leq & \mu(\Um \RmB^{s/2}\phi,\Um\RmB^{s/2} \phi)
         \cdot \mu(\RmB^{s/2}\psi,\RmB^{2(1-s)}\RmB^{s/2} \psi) \\
  & \leq & \mu_s(\phi,\phi)\cdot \mu_s(\psi,\psi)\,, \quad \phi,\psi \in S\,.
\end{eqnarray*}
Here, we have used that $\RmB^{2(1-s)} \leq 1$.
\\[6pt]
(c) If $(\phi_n)$ is a $\mu$-Cauchy-sequence in $H_{\mu}$, then it is,
by continuity of $\RmB^{s/2}$, also a $\mu_s$-Cauchy-sequence in
$H_s$, the $\mu_s$-completion of $S$. Via this identification, we obtain
an embedding $j : H_{\mu} \to H_s$. Notice that $j(\psi) = \psi$
for all $\psi \in S$, so $j$ has dense range; however, one has
\begin{equation}
 \mu_s(j(\phi),j(\psi)) = \mu(\phi,\RmB^s \psi)
\end{equation}
for all $\phi,\psi \in H_{\mu}$. Therefore $j$ need not be injective.
Now let $R_s$ be the polarizator of $\mu_s$. Then we have
\begin{eqnarray*}
2\mu_s(j(\phi),R_s j(\psi))\  = \  \sigma_{\mu}(\phi,\psi) & = &
   2 \mu(\phi,R_{\mu}\psi) \\
& = & 2 \mu(\phi,\RmB^s\Um\RmB^{1-s}\psi) \\
& = & 2 \mu_s(j(\phi),j(\Um\RmB^{1-s})\psi)
         \,,\quad \phi,\psi \in H_{\mu}\,.
\end{eqnarray*}
This yields
\begin{equation}
 R_s \lcrc j = j \lcrc \Um\RmB^{1-s}
\end{equation}
on $H_{\mu}$. Since by assumption $\mu_s$ is pure, we have
$R_s^2 = -1$ on $H_s$, and thus
$$ j = - R_s j \Um\RmB^{1-s} = - j(\Um\RmB^{1-s})^2 \,.$$
By (2.12) we may conclude
$$ \RmB^{2s} = - \Um \RmB \Um \RmB = \Um^*\Um\RmB^2 = \RmB^2 \,, $$
which entails $\RmB^s = \RmB$. Since $\RmB \leq 1$, we see that for
$s \leq r \leq 1$ we have
$$ \RmB = \RmB^s \geq \RmB^r \geq \RmB \,,$$
hence $\RmB^r = \RmB$ for $s \geq r \geq 1$. Whence $\RmB^r = \RmB$
for all $r > 0$. This proves the first part of the statement.

For the second part we observe that $\mu \in \prS$ implies that
$\RmB$, and hence also $\RmB^s$ for $0 < s < 1$, is injective. Then the
equation $\RmB^s = \RmB$ implies that $\RmB^s(\RmB^{1-s} - 1) = 0$,
and by the injectivity of $\RmB^s$ we may conclude $\RmB^{1-s} =1$.
Since $s$ was assumed to be strictly less than 1, it follows that
$\RmB^r = 1$ for all $r \geq 0$; in particular, $\RmB =1$.
\\[6pt]
(d) Assume that $\mu_s$ dominates $\sigma$ for some $s > 1$, i.e.\ it
holds that
$$ 4|\mu(\phi,\Um\RmB\psi)|^2 = |\sigma_{\mu}(\phi,\psi)|^2
\leq 4\cdot \mu(\phi,\RmB^s\phi)\cdot \mu(\psi,\RmB^s\psi)\,, \quad \phi,\psi
\in H_{\mu}\,, $$
which implies, choosing $\phi = \Um \psi$, the estimate
$$ \mu(\psi,\RmB \psi) \leq \mu(\psi,\RmB^s \psi) \,,
 \quad \psi \in H_{\mu}\,,$$
i.e.\ $\RmB \leq \RmB^s$. On the other hand, $\RmB \geq \RmB^r \geq \RmB^s$
holds for all $1 \leq r \leq s$ since $\RmB \leq 1$. This implies
$\RmB^r = \RmB$ for all $r > 0$. For the second part of the statement one
uses the same argument as given in (c). \\[6pt]
(e) In view of (2.13) it holds that
\begin{eqnarray*}
|R_s|^2j & = & - R_s^2 j\  =\  - R_s j \Um \RmB^{1-s} \\
     & = & - j \Um \RmB^{1-s}\Um\RmB^{1-s}\  =\   - j \Um^2 (\RmB^{1-s})^2 \,.
\end{eqnarray*}
Iterating this one has for all $n \in \NN$
$$ |R_s|^{2n} j = (-1)^n j \Um^{2n}(\RmB^{1-s})^{2n}\,. $$
Inserting this into relation (2.12) yields for all $n \in \NN$
\begin{eqnarray}
 \mu_s(j(\phi),|R_s|^{2n}j(\psi)) & = & \mu(\phi,
      \RmB^s (-1)^n \Um^{2n}(\RmB^{1-s})^{2n}\psi) \\
   & = & \mu(\phi,\RmB^s(\RmB^{1-s})^{2n}\psi)\,,\quad \phi,\psi \in H_{\mu}\,.
 \nonumber
\end{eqnarray}
For the last equality we used that $\Um$ commutes with $|R_s|^s$
and $\Um^2\RmB = - \RmB$. Now let $(P_n)$ be a sequence of polynomials
on the intervall $[0,1]$ converging uniformly to the square root
function on $[0,1]$. From (2.14) we infer that
$$ \mu_s(j(\phi),P_n(|R_s|^2)j(\psi)) = \mu(\phi,\RmB^s P_n((\RmB^{1-s})^2)
\psi)\,, \quad \phi, \psi \in H_{\mu} $$
for all $n \in \NN$, which in the limit $n \to \infty$ gives
$$ \mu_s(j(\phi),|R_s|j(\psi)) = \mu(\phi,\RmB\psi)\,, \quad
\phi,\psi \in H_{\mu}\,, $$
as desired. $\Box$
\\[10pt]
 Proposition 2.1 underlines the special role of
$\mT = \mu_1$. Clearly, one has $\mT = \mu$ iff $\mu \in \puS$.
Chmielowski has proved another interesting connection between
$\mu$ and $\mT$ which we briefly mention here. Suppose that
$\{T_t\}$ is a one-parametric group of symplectomorphisms of
$(S,\sigma)$, and let $\{\alpha_t\}$ be the automorphism group
on $\AS$ induced by it via $\alpha_t(W(\phi)) = W(T_t\phi)$,
$\phi \in S,\ t \in \RR$. An $\{\alpha_t\}$-invariant quasifree
state $\om$ on $\AS$ is called {\it regular} if the unitary
group which implements $\{\alpha_t\}$ in the GNS-representation
$(\Hh_{\mu},\pi_{\mu},\Omega_{\mu})$ of $\om$ is strongly
continuous and leaves no non-zero vector in the one-particle space
of $\Hh_{\mu}$ invariant. Here, the one-particle space is spanned
by all vectors of the form $\left. \frac{d}{dt} \right|_{t = 0}
\pi_{\mu}(W(t\phi))\Omega_{\mu}$, $\phi \in S$.
It is proved in [11] that, if $\omega_{\mu}$ is a regular quasifree
KMS-state for $\{\alpha_t\}$, then $\omega_{\mT}$ is the unique
regular quasifree groundstate for $\{\alpha_t\}$. As explained in
[11], the passage from $\mu$ to $\mT$ can be seen as a rigorous
form of ``frequency-splitting'' methods employed in the canonical
quantization of classical fields for which $\mu$ is induced
by the classical energy norm. We shall come back to this in the
concrete example of the Klein-Gordon field in Sec.\ 3.4.

It should be noted that the purification map $\tilde{\cdot} :
\qS \to \puS$, $\mu \mapsto \mT$, assigns to a quasifree state
$\omega_{\mu}$ on $\AS$ the pure quasifree state $\omega_{\mT}$
which is again a state on $\AS$. This is different from the
well-known procedure of assigning to a state $\omega$ on a
$C^*$-algebra $\AA$, whose GNS representation is primary,
 a pure state $\omega_0$ on $\AA^{\circ} \otimes
\AA$.
 ($\AA^{\circ}$ denotes the opposite algebra of $\AA$,
cf.\ [75].) That procedure was introduced by Woronowicz and is an
abstract version of similar constructions for quasifree states on
CCR- or CAR-algebras [45,54,75]. Whether the purification map
$\omega_{\mu} \mapsto \omega_{\mT}$ can be generalized from quasifree states
on CCR-algebras to a procedure of assigning to (a suitable class of)
states on a generic $C^*$-algebra pure states on that same algebra,
is in principle an interesting question, which however we shall not
investigate here.

\begin{Theorem} ${}$\\[6pt]
(a)  Let $H$ be a (real or complex) Hilbertspace with
scalar product $\mu(\,.\,,\,.\,)$, $R$ a (not necessarily bounded) normal
operator in $H$, and $V,W$ two $\mu$-bounded linear operators on $H$
which are $R$-adjoint, i.e.\ they satisfy
\begin{equation}
 W{\rm dom}(R) \subset {\rm dom}(R) \quad {\it and} \quad V^*R = R W \quad
{\rm on \ \   dom}(R) \,.
\end{equation}
Denote by $\mu_s$ the Hermitean form on ${\rm dom}(|R|^{s/2})$ given by
$$ \mu_s(x,y) := \mu(|R|^{s/2}x,|R|^{s/2} y)\,, \quad
 x,y \in {\rm dom}(|R|^{s/2}),\ 0 \leq s \leq 2\,.$$ We write
$||\,.\,||_0 := ||\,.\,||_{\mu} := \mu(\,.\,,\,.\,)^{1/2}$ and
$||\,.\,||_s := \mu_s(\,.\,,\,.\,)^{1/2}$ for the corresponding semi-norms.

Then it holds for all $0 \leq s \leq 2$ that
$$ V{\rm dom}(|R|^{s/2}) \subset {\rm dom}(|R|^{s/2}) \quad
{\it and} \quad W{\rm dom}(|R|^{s/2}) \subset {\rm dom}(|R|^{s/2}) \,, $$
and $V$ and $W$ are $\mu_s$-bounded for $0 \leq s \leq 2$.
More precisely, the estimates
\begin{equation}
||\,Vx\,||_0 \leq v\,||\,x\,||_0 \quad {\rm and} \quad
||\,Wx\,||_0 \leq w\,||\,x\,||_0\,, \quad x \in H\,,
\end{equation}
with suitable constants $v,w > 0$, imply that
\begin{equation}
||\,Vx\,||_s \leq w^{s/2}v^{1 -s/2}\,||\,x\,||_s \quad {\rm and} \quad
||\,Wx\,||_s \leq v^{s/2}w^{1-s/2}\,||\,x\,||_s \,,
\end{equation}
for all
$  x \in {\rm dom}(|R|^{s/2})$ and
 $0 \leq s \leq 2$.
\\[6pt]
(b)\ \ \   (Corollary of (a))\ \ \ \ 
Let $(S,\sigma)$ be a symplectic space, $\mu \in \qS$ a dominating
scalar product on $(S,\sigma)$, and $\mu_s$, $0 < s \leq 2$, the
scalar products on $S$ defined in (2.11). Then we have relative
$\mu-\mu_s$ continuity of each pair 
$V,W$ of symplectically adjoint linear maps of $(S,\sigma)$
for all $0 < s \leq 2$. More precisely, for each pair $V,W$ of
symplectically adjoint linear maps of $(S,\sigma)$, the estimates
(2.16) for all $x \in S$ imply (2.17) for all $x \in S$.
\end{Theorem}
{\it Remark.} (i) In view of the fact that the operator  $R$
of part (a) of the Theorem may be unbounded, part (b) can  be
extended to situations where it is not assumed that the scalar
product $\mu$ on $S$ dominates the symplectic form $\sigma$.
\\[6pt]
(ii) When it is additionally assumed that $V = T$ and $W = T^{-1}$
with symplectomorphisms $T$ of $(S,\sigma)$, we refer in that
case to the situation of relative continuity of the pairs
$V,W$ as relative continuity of symplectomorphisms.
 In Example 2.3 after the proof of Thm.\ 2.2 we show that
relative $\mT - \mu$ continuity of symplectomorphisms fails in general.
Also, it is not the case that relative $\mu - \mu'$ continuity
of symplectomorphisms holds if $\mu'$ is an arbitrary element
in $\puS$ which is dominated by $\mu$ ($||\,\phi\,||_{\mu'} \leq
{\rm const.}||\,\phi\,||_{\mu}$, $\phi \in S$), see Example 2.4
below. This shows that the special relation between $\mu$ and $\mT$
(resp., $\mu$ and the $\mu_s$) expressed in (2.11,2.15) is important
for the derivation of the Theorem.
\\[10pt]
{\it Proof of Theorem 2.2.} (a)\ \ \  In a first step, let it be 
supposed that $R$ is bounded.
 From the assumed relation (2.15)
and its adjoint relation $R^*V = W^* R^*$ we obtain, for $\epsilon' > 0$
arbitrarily chosen,
\begin{eqnarray*}
V^* (|R|^2 + \epsilon' 1) V & = & V^*RR^* V + \epsilon' V^*V \ =
              \ RWW^*R^* + \epsilon' V^*V \\
& \leq & w^2 |R|^2 + \epsilon' v^21 \  \leq \  w^2(|R|^2 + \epsilon 1)
\end{eqnarray*}
with $\epsilon := \epsilon'v^2/w^2$.
 This entails for the operator norms
$$ ||\,(|R|^2 + \epsilon'  1 )^{1/2}V \,||
\  \leq\  w\,||\,(|R|^2 + \epsilon 1)^{1/2}\,|| \,, $$
and since $(|R|^2 + \epsilon  1)^{1/2}$ has a bounded inverse,
$$ ||\,(|R|^2 + \epsilon'  1)^{1/2} V
 (|R|^2 + \epsilon  1 )^{-1/2} \,||\  \leq\  w\,. $$
On the other hand, one clearly has
$$ ||\,(|R|^2 + \epsilon'  1)^0 V (|R|^2 + \epsilon  1)^0\,||
\  =\  ||\,V\,||\  \leq\   v\,. $$
Now these estimates are preserved if $R$ and $V$
are replaced by their complexified versions on the complexified
Hilbertspace $H \oplus iH = \CC \otimes H$.
Thus, identifying if necessary
 $R$ and $V$ with their complexifications, a standard interpolation
argument (see Appendix A) can be applied to yield
$$ ||\,(|R|^2 + \epsilon'  1)^{\alpha} V
      (|R|^2 + \epsilon  1)^{-\alpha} \,||\  \leq\  w^{2\alpha}
v^{1 - 2\alpha}  $$
for all $0 \leq \alpha \leq 1/2$. Notice that this inequality
holds uniformly in $\epsilon' > 0$. Therefore we may conclude that
$$ ||\,|R|^{2\alpha}Vx \,||_0\  \leq\  w^{2\alpha}v^{1 - 2\alpha}
\,||\,|R|^{2\alpha}x\,||_0
\,, \quad x \in H\,,\ 0 \leq \alpha \leq 1/2\,,$$
which is the required estimate for $V$. The analogous bound for
$W$ is obtained through replacing $V$ by $W$ in the
given arguments. 

Now we have to extend the argument to the case that $R$ is unbounded.
Without restriction of generality we may assume that the Hilbertspace
$H$ is complex, otherwise we complexify it and with it all the
operators $R$,$V$,$W$, as above, thereby preserving their assumed properties.
Then let $E$ be the spectral measure of $R$, and denote by
$R_r$ the operator $E(B_r)RE(B_r)$ where $B_r := \{z \in \CC:
|z| \leq r\}$, $r > 0$. Similarly define $V_r$ and $W_r$. From the
assumptions it is seen that $V_r^*R_r = R_rW_r$ holds for all
$r >0$. Applying the reasoning of the first step we arrive, for each
$0 \leq s \leq 2$, at the bounds
$$ ||\,V_r x\,||_s \leq w^{s/2}v^{1-s/2}\,||\,x\,||_s \quad {\rm and}
\quad ||\,W_r x\,||_s \leq v^{s/2}w^{1-s/2}\,||\,x\,||_s \,,$$
which hold uniformly in $r >0$ for all $x \in {\rm dom}(|R|^{s/2})$.
From this, the statement of the Proposition follows.\\[6pt]
(b) This is just an application of (a), identifying  $H_{\mu}$ with $H$,
$R_{\mu}$ with $R$ and $V,W$ with their bounded extensions to $H_{\mu}$.
$\Box$ 
\\[10pt]
{\bf Example 2.3} We exhibit a symplectic space $(S,\sigma)$
with $\mu \in \prS$ and a symplectomorphism $T$ of $(S,\sigma)$
where $T$ and $T^{-1}$ are continuous with respect to $\mT$,
but not with respect to $\mu$. \\[6pt]
Let $S := {\cal S}(\RR,\CC)$, the Schwartz space of rapidly decreasing
testfunctions on $\RR$, viewed as real-linear space.
By $\langle \phi,\psi \rangle := \int \overline{\phi} \psi \,dx$
we denote the standard $L^2$ scalar product. As a symplectic form on
$S$ we choose
$$ \sigma(\phi,\psi) := 2 {\sf Im}\langle \phi,\psi \rangle\,, \quad
\phi,\psi \in S\,. $$
Now define on $S$ the strictly positive, essentially selfadjoint operator
$A\phi := -\frac{d^2}{dx^2}\phi + \phi$,
 $\phi \in S$, in $L^2(\RR)$. Its closure
will again be denoted by $A$; it is bounded below by $1$.
A real-linear scalar product $\mu$ will be defined on $S$ by
$$ \mu(\phi,\psi) := {\sf Re}\langle A\phi,\psi \rangle\,, \quad \phi,\psi \in
S. $$
Since $A$ has lower bound $1$, clearly $\mu$ dominates $\sigma$, and
one easily obtains $R_{\mu} = - i A^{-1}$, $\RmB = A^{-1}$.
Hence $\mu \in \prS$ and
$$ \mT(\phi,\psi) = {\sf Re}\langle \phi,\psi \rangle\,,
           \quad   \phi,\psi \in S\,.$$
Now consider the operator
$$ T : S \to S\,, \quad \ \  (T\phi)(x) := {\rm e}^{-ix^2}\phi(x)\,,
\ \ \  x \in \RR,\  \phi \in S\,. $$
Obviously $T$ leaves the $L^2$ scalar product invariant, and hence also
$\sigma$ and $\mT$. The inverse of $T$ is just $(T^{-1}\phi)(x) =
{\rm e}^{i x^2}\phi(x)$, which of course leaves $\sigma$ and $\mT$
invariant as well. However, $T$ is not continuous with respect to $\mu$.
To see this, let $\phi \in S$ be some non-vanishing smooth function
with compact support, and define
$$ \phi_n(x) := \phi(x -n)\,, \quad x \in \RR, \   n \in \NN\,. $$
Then $\mu(\phi_n,\phi_n) = {\rm const.} > 0$ for all $n \in \NN$.
We will show that $\mu(T\phi_n,T\phi_n)$ diverges for $n \to \infty$.
We have
\begin{eqnarray}
 \mu(T\phi_n,T\phi_n) & = & \langle A T\phi_n,T\phi_n \rangle
       \geq \int \overline{(T\phi_n)'}(T\phi_n)' \, dx   \\
       & \geq & \int (4x^2|\phi_n(x)|^2 + |\phi_n'(x)|^2)\, dx
        - \int 4 |x \phi_n'(x)\phi_n(x)|\,dx \,,\nonumber
\end{eqnarray}
where the primes indicate derivatives and we have used that
$$ |(T\phi_n)'(x)|^2 = 4x^2|\phi_n(x)|^2 + |\phi_n'(x)|^2
                        + 4\cdot {\sf Im}(ix \overline{\phi_n}(x)\phi_n'
  (x))\,. $$
Using
 a substitution of variables, one can see that in the last term
of (2.18) the positive integral grows like $n^2$ for large $n$, thus
dominating eventually the negative integral which grows only like $n$.
So $\mu(T\phi_n,T\phi_n) \to \infty$ for $n \to \infty$, showing that
$T$ is not $\mu$-bounded.
\\[10pt]
{\bf Example 2.4} We give an example of a symplectic space
$(S,\sigma)$, a $\mu \in \prS$ and a $\mu' \in \puS$, where
$\mu$ dominates $\mu'$ and where there is a symplectomorphism $T$
of $(S,\sigma)$ which together with its inverse is $\mu$-bounded,
but not $\mu'$-bounded.\\[6pt]
We take $(S,\sigma)$ as in the previous example and write for each
$\phi \in S$, $\phi_0 := {\sf Re}\phi$ and $\phi_1 := {\sf Im}\phi$.
The real scalar product $\mu$ will be defined by
$$ \mu(\phi,\psi) := \langle\phi_0,A\psi_0\rangle + \langle \phi_1,
\psi_1 \rangle \,, \quad \phi,\psi \in S\,, $$
where the operator $A$ is the same as in the example before. Since its
lower bound is $1$, $\mu$ dominates $\sigma$, and it is not difficult
to see that $\mu$ is even primary. The real-linear scalar product
$\mu'$ will be taken to be
$$ \mu'(\phi,\psi) = {\sf Re}\langle \phi,\psi \rangle\,, \quad
\phi,\psi \in S\,.$$
We know from the example above that $\mu' \in \puS$. Also, it is
clear that $\mu'$ is dominated by $\mu$. Now consider the
real-linear map $T: S \to S$ given by
$$ T(\phi_0 + i\phi_1) := A^{-1/2} \phi_1 - i A^{1/2}\phi_0\,, \quad
\phi \in S\,.$$
One checks easily that this map is bijective with $T^{-1} = - T$,
and that $T$ preserves the symplectic form $\sigma$. Also,
$||\,.\,||_{\mu}$ is preserved by $T$ since
$$ \mu(T\phi,T\phi) = \langle \phi_1,\phi_1 \rangle +
              \langle A^{1/2}\phi_0,A^{1/2}\phi_0 \rangle = \mu(\phi,\phi)\,,
\quad \phi \in S\,.$$
On the other hand, we have for each $\phi \in S$
$$ \mu'(T\phi,T\phi) = \langle \phi_1,A\phi_1 \rangle + \langle \phi_0,
A^{-1}\phi_0 \rangle \,, $$
and this expression is not bounded by a ($\phi$-independent) constant times
$\mu'(\phi,\phi)$, since $A$ is unbounded with respect to the
$L^2$-norm.
\newpage
\section{The Algebraic Structure of Hadamard Vacuum Representations}
\setcounter{equation}{0}
${}$
\\[20pt]
{\bf 3.1 Summary of Notions from Spacetime-Geometry}
\\[16pt]
We recall that a spacetime manifold consists of a pair
$(M,g)$, where $M$ is a smooth, paracompact,  four-dimensional
manifold without boundaries, and $g$ is a Lorentzian metric for $M$
with signature $(+ - -\, - )$. (Cf.\ [33,52,70], see these references
also for further discussion of the notions to follow.)
It will be assumed that $(M,g)$ is time-orientable, and moreover,
globally hyperbolic. The latter means that $(M,g)$ possesses
Cauchy-surfaces, where by a Cauchy-surface
we always mean  a {\it smooth}, spacelike
hypersurface which is intersected exactly once by each inextendable
causal curve in $M$. It can be shown [15,28] that this is
equivalent to the statement that $M$ can be smoothly foliated in
Cauchy-surfaces. Here, a foliation of $M$ in Cauchy-surfaces is
a diffeomorphism $F: \RR \times \Sigma \to M$, where $\Sigma$ is a
smooth 3-manifold so that $F(\{t\} \times \Sigma)$ is, for each
$t \in \RR$, a Cauchy-surface, and the curves 
$t \mapsto F(t,q)$
 are timelike for all $q \in\Sigma$.
(One can even show that, if global hyperbolicity had been
defined by requiring only the existence of a non necessarily
smooth or spacelike Cauchy-surface (i.e.\ a topological
hypersurface which is intersected exactly once by each 
inextendable causal curve), then it is still true that a globally
hyperbolic spacetime can be smoothly foliated in Cauchy-surfaces,
see [15,28].)

 We shall also be interested in
ultrastatic globally hyperbolic spacetimes.
 A globally hyperbolic spacetime
is said to be {\it ultrastatic} if a foliation
$F : \RR \times \Sigma \to M$ in Cauchy-surfaces can be found so that
$F_*g$ has the form $dt^2 \oplus (- \gamma)$ with a complete
($t$-independent) Riemannian metric $\gamma$ on $\Sigma$.
This particular foliation will then be called a {\it natural foliation}
of the ultrastatic spacetime. (An ultrastatic spacetime may posses
more than one natural foliation, think e.g.\ of Minkowski-spacetime.)

The notation for the causal sets and domains of dependence
will be recalled: Given a spacetime $(M,g)$ and
$\Oo \subset M$, the set $J^{\pm}(\Oo)$ (causal future/past of
$\Oo$) consists of all points $p \in M$ which can be reached by
future/past directed causal curves emanating from $\Oo$. The set
$D^{\pm}(\Oo)$ (future/past domain of dependence of $\Oo$) is defined
as consisting of all $p \in J^{\pm}(\Oo)$ such that every
past/future inextendible causal curve starting at $p$ intersects $\Oo$.
One writes $J(\Oo) := J^+(\Oo) \cup J^-(\Oo)$ and $D(\Oo) :=
D^+(\Oo) \cup D^-(\Oo)$. They are called the {\it causal set}, and
the {\it domain of dependence}, respectively, of $\Oo$.

For $\Oo \subset M$, we denote by $\Oo^{\perp} := {\rm int}(M
\backslash J(\Oo))$ the {\it causal complement} of $\Oo$,
 i.e.\ the largest {\it open}
set of points which cannot be connected to $\Oo$ by any causal curve.

A set of the form $\Oo_G := {\rm int}\,D(G)$, where $G$ is a subset
of some Cauchy-surface $\Sigma$ in $(M,g)$, will be referred to
as the {\it diamond based on} $G$; we shall also say that
$G$ is the {\it base} of $\Oo_G$. We note that if $\Oo_G$ is a
diamond, then $\Oo_G^{\perp}$ is again a diamond, based on
$\Sigma \backslash \overline{G}$.
 A diamond will be called
{\it regular} if $G$ is an open, relatively compact subset of
$\Sigma$ and if the boundary $\partial G$ of $G$ is contained in
the union of finitely many smooth, two-dimensional submanifolds
of $\Sigma$. 

Following [45], we say that an open neighbourhood $N$ of a
Cauchy-surface $\Sigma$ in $(M,g)$ is a {\it causal normal neighbourhood}
of $\Sigma$ if (1) $\Sigma$ is a Cauchy-surface for $N$, and
(2) for each pair of points $p,q \in N$ with $p \in J^+(q)$, there
is a convex normal neighbourhood $\Oo \subset M$ such that 
$J^-(p) \cap J^+(q) \subset \Oo$. Lemma 2.2 of [45] asserts the
existence of causal normal neighbourhoods for any Cauchy-surface
$\Sigma$.
\\[20pt]
{\bf 3.2 Some Structural Aspects of
Quantum Field Theory in Curved Spacetime}
\\[16pt]
In the present subsection, we shall address some of the problems
one faces in the formulation of quantum field theory in curved
spacetime, and explain the notions of local definiteness, local
primarity, and Haag-duality. In doing so, we follow our presentation
in [67] quite closely. Standard general references related to the
subsequent discussion are [26,31,45,71].

Quantum field theory in curved spacetime (QFT in CST, for short)
means that one considers quantum field theory means that one considers
quantum fields propagating in a (classical) curved background spacetime
manifold $(M,g)$. In general, such a spacetime need not possess
any symmetries, and so one cannot tie the notion of ``particles''
or ``vacuum'' to spacetime symmetries, as one does in quantum field
theory in Minkowski spacetime. Therefore, the problem of
how to characterize the physical states arises. For the discussion
of this problem, the setting of algebraic quantum field theory is
particularly well suited. Let us thus summarize some of the relevant
concepts of algebraic QFT in CST. Let a spacetime manifold
$(M,g)$ be given. The observables of a quantum system (e.g.\ a quantum
field) situated in $(M,g)$ then have the basic structure of a map
$\Oo \to \AO$, which assigns to each open, relatively compact
subset  $\Oo$ of $M$ a $C^*$-algebra $\AO$,\footnote{
Throughout the paper, $C^*$-algebras are assumed to be unital, i.e.\
to possess a unit element, denoted by ${ 1}$. It is further
assumed that the unit element is the same for all the $\AO$.}
 with the properties:\footnote{where $[\AA(\Oo_1),\AA(\Oo_2)]
= \{A_1A_2 - A_2A_1 : A_j \in \AA(\Oo_j),\ j =1,2 \}$.}
\begin{equation}
{\it Isotony:}\quad \quad \Oo_1 \subset \Oo_2
\Rightarrow \AA(\Oo_1) \subset \AA(\Oo_2)
\end{equation}
\begin{equation}
{\it Locality:} \quad \quad \Oo_1 \subset \Oo_2^{\perp}
\Rightarrow [\AA(\Oo_1),\AA(\Oo_2)] = \{0 \} \,.
\end{equation}
A map $\Oo \to \AO$ having these properties is called a {\it net
of local observable algebras} over $(M,g)$. We recall that the
conditions of locality and isotony are motivated by the idea
that each $\AO$ is the $C^*$-algebra formed by the observables
which can be measured within the spacetime region $\Oo$ on the
system. We refer to [31] and references given there for further
discussion.

The collection of all open, relatively compact subsets of
$M$ is directed with respect to set-inclusion, and so we can, in view
of (3.1), form the smallest $C^*$-algebra $\AA := \overline{
\bigcup_{\Oo}\AO}^{||\,.\,||}$ which contains all local algebras
$\AO$.
For the description of a system we need not only observables but also
states. The set $\AA^{*+}_1$ of all positive, normalized linear
functionals on $\AA$ is mathematically referred to as the set
of {\it states} on $\AA$, but not all elements of $\AA^{*+}_1$
represent physically realizable states of the system. Therefore, given
a local net of observable algebras $\Oo \to \AO$ for a physical
system over $(M,g)$, one must specify the set of physically relevant states
$\Ss$, which is a suitable subset of $\AA^{*+}_1$.

We have already mentioned in Chapter 2 that every state $\omega \in
\AA^{*+}_1$ determines canonically its GNS representation
$(\Hh_{\omega},\pi_{\omega},\Omega_{\omega})$ and thereby induces
a net of von Neumann algebras (operator algebras on $\Hh_{\omega}$)
$$ \Oo \to \Rr_{\omega}(\Oo) := \pi_{\omega}(\Oo)^- \,. $$
Some of the mathematical properties of the GNS representations, and of
the induced nets of von Neumann algebras, of states $\omega$ on
$\AA$ can naturally be interpreted physically. Thus one obtains
constraints on the states $\omega$ which are to be viewed as
physical states. Following this line of thought, Haag, Narnhofer
and Stein [32] formulated what they called the ``principle of local
definiteness'', consisting of the following three conditions
to be obeyed by any collection $\Ss$ of physical states.
\\[10pt]
{\bf Local Definiteness:} ${}\ \ \bigcap_{\Oo \owns p}
\Rr_{\omega}(\Oo) = \CC \cdot { 1}$ for all $\omega \in \Ss$
and all $p \in M$.
\\[6pt]
{\bf Local Primarity:} \ \ For each $\omega \in \Ss$, $\Rr_{\omega}
(\Oo)$ is a factor.
\\[6pt]
{\bf Local Quasiequivalence:} For each pair $\omega_1,\omega_2 \in \Ss$
and each relatively compact, open $\Oo \subset M$, the representations
$\pi_{\omega_1} | \AO$ and $\pi_{\omega_2} | \AO$ of $\AO$ are
quasiequivalent.
\\[10pt]
{\it Remarks.} (i) We recall  (cf.\ the first Remark in Section 2) that
$\Rr_{\omega}(\Oo)$ is a factor if $\Rr_{\omega}(\Oo) \cap \Rr_{\omega}
(\Oo)' = \CC \cdot { 1}$ where the prime means taking the
commutant. We have not stated in the formulation of local primarity
for which regions $\Oo$ the algebra $\Rr_{\omega}(\Oo)$ is required to be
a factor. The regions $\Oo$ should be taken from a class of subsets of
$M$ which forms a base for the topology.
\\[6pt]
(ii) Quasiequivalence of representations means unitary equivalence up to
multiplicity. Another characterization of quasiequivalence is to say
that the folia of the representations coincide, where the
{\it folium} of
a representation $\pi$ is defined as the set of all $\omega \in
\AA^{*+}_1$ which can be represented as $\omega(A) = tr(\rho\,\pi(A))$
with a density matrix $\rho$ on the representation Hilbertspace of $\pi$.
\\[6pt]
(iii) Local definiteness and quasiequivalence together express
that physical states have finite (spatio-temporal) energy-density
with respect to each other, and local primarity and quasiequivalence
rule out local macroscopic observables and local superselection
rules. We refer to [31] for further discussion and background
material.

A further, important property which one expects to be satisfied
for physical states $\omega \in \Ss$ whose GNS representations are
irreducible \footnote{It is easy to see that,
in the presence of local primarity, Haag-duality will be
violated if $\pi_{\omega}$ is not irreducible.} is
\\[10pt]
{\bf Haag-Duality:} \ \  $\Rr_{\omega}(\Oo^{\perp})' = \Rr_{\omega}(\Oo)$, \\
which should hold for the causally complete regions $\Oo$, i.e.\ those
satisfying $(\Oo^{\perp})^{\perp} = \Oo$, where $\Rr_{\omega}(\Oo^{\perp})$
is defined as the von Neumann algebra generated by all the $\Rr_{\omega}
(\Oo_1)$ so that $\overline{\Oo_1} \subset \Oo^{\perp}$.
\\[10pt]
We comment that Haag-duality means that the von Neumann algebra
$\Rr_{\omega}(\Oo)$ of local observables is maximal in the sense
that no further observables can be added without violating the
condition of locality. It is worth mentioning here that the condition
of Haag-duality plays an important role in the theory of superselection
sectors in algebraic quantum field theory in Minkowski spacetime
[31,59]. For local nets of observables generated by Wightman fields
on Minkowski spacetime it follows from the results of Bisognano and
Wichmann [4] that a weaker condition of ``wedge-duality'' is
always fulfilled, which allows one to pass to a new, potentially
larger local net (the ``dual net'') which satisfies Haag-duality.

In quantum field theory in Minkowski-spacetime where one is given
a vacuum state $\omega_0$, one can define the set of physical states
$\Ss$ simply as the set of all states on $\AA$ which are locally
quasiequivalent (i.e., the GNS representations of the states are
locally quasiequivalent to the vacuum-representation) to $\omega_0$.
It is obvious that local quasiequivalence then holds for $\Ss$.
Also, local definiteness holds in this case, as was proved by Wightman
[72]. If Haag-duality holds in the vacuum representation (which,
as indicated above, can be assumed to hold quite generally), then it
does not follow automatically that all pure states locally quasiequivalent
to $\omega_0$ will also have GNS representations fulfilling Haag-duality;
 however, it follows once some regularity conditions are satisfied
which have been checked in certain quantum field models [19,61].
So far there seems to be no general physically motivated criterion
enforcing local primarity of a quantum field theory in algebraic
formulation in Minkowski spacetime. But it is known that many quantum
field theoretical models satisfy local primarity.

For QFT in CST we do in general not know what a vacuum state is and so
$\Ss$ cannot be defined in the same way as just described. Yet in some
cases (for some quantum field models) there may be a set
$\Ss_0 \subset \AA^{*+}_1$ of distinguished states, and if this class
of states satisfies the four conditions listed above, then the set
$\Ss$, defined as consisting of all states $\omega_1 \in \AA^{*+}_1$
which are locally quasiequivalent to any (and hence all) $\omega
\in \Ss_0$, is a good candidate for the set of physical states.

For the free scalar Klein-Gordon field (KG-field) on
a globally hyperbolic spacetime, the following classes of
states have been suggested as distinguished, physically
reasonable states \footnote{The following list is not meant
to be complete, it comprises some prominent families of states
of the KG-field over a generic class of spacetimes for which
mathematically sound results are known. Likewise, the indicated
references are by no means exhaustive.}
\begin{itemize}
\item[(1)] (quasifree) states fulfilling local stability
[3,22,31,32]
\item[(2)] (quasifree) states fulfilling the wave front set (or microlocal)
spectrum condition [6,47,55]
\item[(3)] quasifree Hadamard states [12,68,45]
\item[(4)] adiabatic vacua [38,48,53]
\end{itemize}
The list is ordered in such a way that the less restrictive
condition preceeds the stronger one. There are a couple of
comments to be made here. First of all, the specifications
(3) and (4) make use of the information that one deals with
the KG-field (or at any rate, a free field obeying a linear
equation of motion of hyperbolic character), while the
conditions  (1) and (2) do not require such input and are
applicable to general -- possibly interacting -- quantum
fields over curved spacetimes.
(It should however be mentioned that only for the KG-field (2)
is known to be stronger than (1). The relation between (1) and (2)
for more general theories is not settled.)
 The conditions imposed on the classes of
states (1), (2) and (3) are related in that they are ultralocal remnants
of the spectrum condition requiring a certain regularity of the
short distance behaviour of the respective states which can be
formulated in generic spacetimes. The class of states (4) is
more special and can only be defined for the KG-field (or other
linear fields) propagating in Robertson-Walker-type spacetimes.
Here a distinguished choice of a time-variable can be made,
and the restriction imposed on adiabatic vacua is a regularity
condition on their spectral behaviour with respect to that
special choice of time. (A somewhat stronger formulation of
local stability has been proposed in [34].) 

It has been found by Radzikowski [55] that for quasifree states of the
KG-field over generic globally hyperbolic spacetimes
the classes (2) and (3) coincide. The microlocal spectrum condition
is further refined and applied in [6,47]. Recently it was
proved by Junker [38] that adiabatic vacua of the KG-field
in Robertson-Walker spacetimes fulfill the microlocal spectrum
condition and thus are, in fact, quasifree Hadamard states.
The notion of the microlocal spectrum condition and the just mentioned
results related to it draw on pseudodifferential operator
techniques, particularly the notion of the wave front set, see [20,36,37].

 Quasifree Hadamard states of the KG-field (see definition
in Sec.\ 3.4  below)
have been investigated for quite some time. One of the early
studies of these states is [12]. The importance of these
states, especially in the context of the semiclassical
Einstein equation, is stressed in [68]. Other significant
references include [24,25] and, in particular, [45] where,
apparently for the first time, a satisfactory definition of
the notion of a globally Hadamard state is given, cf.\
 Section 3.4 for more details. In [66] it is proved that the
class of quasifree Hadamard states of the KG-field fulfills
local quasiequivalence in generic globally hyperbolic spacetimes
and local definiteness, local primarity and Haag-duality
for the case of ultrastatic globally hyperbolic spacetimes.
As was outlined in the beginning, the purpose of the present
chapter is to obtain these latter results also for arbitrary
globally hyperbolic spacetimes which are not necessarily
ultrastatic. It turns out that some of our previous results
can be sharpened, e.g.\ the local quasiequivalence specializes
in most cases to local unitary equivalence, cf.\ Thm.\ 3.6.
 For a couple of other
results about the algebraic structure of the KG-field as well as other
fields over curved spacetimes we refer to
 [2,6,15,16,17,40,41,46,63,64,65,66,74].
\\[24pt]
{\bf 3.3 The Klein-Gordon Field}
\\[18pt]
In the present section we summarize the quantization of the
classical KG-field over a globally hyperbolic spacetime in the
$C^*$-algebraic formalism. This follows in major parts the
the work of Dimock [16], cf.\ also references given there.
Let $(M,g)$ be a globally hyperbolic spacetime. The KG-equation
with potential term $r$ is
\begin{equation}
(\nabla^a \nabla_a + r) \varphi = 0
\end{equation}
where $\nabla$ is the Levi-Civita derivative of the metric $g$,
the potential function
$r \in C^{\infty}(M,\RR)$ is arbitrary but fixed, and the
sought for solutions $\varphi$ are smooth and real-valued.
Making use of the fact that $(M,g)$ is globally hyperbolic
and drawing on earlier results by Leray, it is shown in
[16] that there are two uniquely determined, continuous
\footnote{ With respect to the usual locally convex topologies
on $C_0^{\infty}(M,\RR)$ and $C^{\infty}(M,\RR)$, cf.\
[13].}
linear maps $E^{\pm}: C_0^{\infty}(M,\RR) \to C^{\infty}(M,\RR)$
with the properties
$$ (\nabla^a \nabla_a + r)E^{\pm}f = f = E^{\pm}(\nabla^a
\nabla_a + r)f\,,\quad f \in C_0^{\infty}(M,\RR)\,, $$
and
$$ {\rm supp}(E^{\pm}f) \subset J^{\pm}({\rm supp}(f))\,,\quad
f \in C_0^{\infty}(M,\RR)\,. $$
The maps $E^{\pm}$ are called the advanced(+)/retarded(--)
fundamental solutions of the KG-equation with potential term
$r$ in $(M,g)$, and their difference $E := E^+ - E^-$ is referred
to as the {\it propagator} of the KG-equation.

One can moreover show that the Cauchy-problem for the KG-equation
is well-posed. That is to say, if $\Sigma$ is any Cauchy-surface
in $(M,g)$, and
$u_0 \oplus u_1 \in C_0^{\infty}(M,\RR) \oplus
C_0^{\infty}(M,\RR)$ is any pair of Cauchy-data on $\Sigma$,
then there exists precisely one smooth solution $\varphi$
of the KG-equation (3.3) having the property that
\begin{equation}
P_{\Sigma}(\varphi) := \varphi | \Sigma \oplus
n^a \nabla_a \varphi|\Sigma = u_0 \oplus u_1\,.
\end{equation}
The vectorfield $n^a$ in (3.4) is the future-pointing
unit normalfield of $\Sigma$.
Furthermore, one has ``finite propagation speed'', i.e.\
 when the supports of $u_0$ and $u_1$ are contained in a subset
$G$ of $\Sigma$, then ${\rm supp}(\varphi) \subset J(G)$. Notice that
compactness of $G$ implies that $J(G) \cap \Sigma'$ is compact
for any Cauchy-surface $\Sigma'$.

 The well-posedness of the
Cauchy-problem is a consequence of the classical energy-estimate
for solutions of second order hyperbolic partial differential
equations, cf.\ e.g.\ [33]. To formulate it, we introduce further
notation. Let $\Sigma$ be a Cauchy-surface for $(M,g)$, and
$\gamma_{\Sigma}$ the Riemannian metric, induced by the ambient
Lorentzian metric, on $\Sigma$. Then denote the Laplacian
operator on $C_0^{\infty}(\Sigma,\RR)$ corresponding to
$\gamma_{\Sigma}$ by $\Delta_{\gamma_{\Sigma}}$, and define the
{\it classical energy scalar product} on
$C_0^{\infty}(\Sigma,\RR) \oplus C_0^{\infty}(\Sigma,\RR)$ by
\begin{equation}
\mu_{\Sigma}^E(u_0 \oplus u_1,
v_0 \oplus v_1 ) := \int_{\Sigma} (u_0 (- \Delta_{\gamma_{\Sigma}} + 1)v_0
+ u_1 v_1) \, d\eta_{\Sigma} \,,
\end{equation}
where $d\eta_{\Sigma}$ is the metric-induced volume measure on
$\Sigma$. As a special case of the energy estimate
presented in [33] one then obtains
\begin{Lemma}
(Classical energy estimate for the KG-field.)
Let $\Sigma_1$ and $\Sigma_2$ be a pair of Cauchy-surfaces
in $(M,g)$ and $G$ a compact subset of $\Sigma_1$. Then there
are two positive constants $c_1$ and $c_2$ so that there
holds the estimate
\begin{equation}
c_1\,\mu^E_{\Sigma_1}(P_{\Sigma_1}(\varphi),P_{\Sigma_1}
(\varphi)) \leq \mu^E_{\Sigma_2}(P_{\Sigma_2}(\varphi),
P_{\Sigma_2}(\varphi)) \leq
c_2 \, \mu^E_{\Sigma_1}(P_{\Sigma_1}(\varphi),P_{\Sigma_1}
(\varphi))
\end{equation}
for all solutions $\varphi$ of the KG-equation (3.3)
which have the property that the supports of the Cauchy-data
$P_{\Sigma_1}(\varphi)$ are contained in $G$.
\footnote{ {\rm The formulation given here is to some extend
more general than the one appearing in [33] where it is 
assumed that $\Sigma_1$ and $\Sigma_2$ are members of a foliation.
However, the more general formulation can be reduced to that case.}}
\end{Lemma}
We shall now indicate that the space of smooth solutions of the
KG-equation (3.3) has the structure of a symplectic space,
locally as well as globally, which comes in several equivalent
versions. To be more specific, observe first that the
Cauchy-data space
$$ \DS := \CoS \oplus \CoS $$
of an arbitrary given Cauchy-surface $\Sigma$ in $(M,g)$
carries a symplectic form
$$ \delta_{\Sigma}(u_0 \oplus u_1, v_0 \oplus v_1)
:= \int_{\Sigma}(u_0v_1 - v_0u_1)\,d\eta_{\Sigma}\,.$$
It will also be observed that this symplectic form
is dominated by the classical energy scalar product
$\mu^E_{\Sigma}$.

                Another symplectic space is $S$, the set
of all real-valued $C^{\infty}$-solutions $\varphi$ of the
KG-equation (3.3) with the property that, given any Cauchy-surface
$\Sigma$ in $(M,g)$, their Cauchy-data $P_{\Sigma}(\varphi)$
have compact support on $\Sigma$. The symplectic form on
$S$ is given by
$$ \sigma(\varphi,\psi) := \int_{\Sigma}(\varphi n^a \nabla_a\psi
-\psi n^a \nabla_a \varphi)\,d\eta_{\Sigma} $$
which is independent of the choice of the Cauchy-surface $\Sigma$
on the right hand side over which the integral is formed;
$n^a$ is again the future-pointing unit normalfield of $\Sigma$.
One clearly finds that for each Cauchy-surface $\Sigma$ the
map $P_{\Sigma} : S \to \DS$ establishes a symplectomorphism
between the symplectic spaces $(S,\sigma)$ and $(\DS,\delta_{\Sigma})$.

A third symplectic space equivalent to the previous ones is obtained
as the quotient $K := \CoM /{\rm ker}(E)$ with symplectic form
$$ \kappa([f],[h]) := \int_M f(Eh)\,d\eta \,, \quad
f,h \in \CoM\,, $$
where $[\,.\,]$ is the quotient map $\CoM \to K$ and
$d\eta$ is the metric-induced volume measure on $M$. 

Then define for any open subset $\Oo \subset M$ with compact
closure the set $K(\Oo) := [C_0^{\infty}(\Oo,\RR)]$. 
One can see that the space $K$ has naturally the structure of
an isotonous, local net $\Oo \to K(\Oo)$ of subspaces, where
locality means that the symplectic form $\kappa([f],[h])$
vanishes for $[f] \in K(\Oo)$ and $[h] \in K(\Oo_1)$
whenever $\Oo_1 \subset \Oo^{\perp}$.
Dimock has proved in [16 (Lemma A.3)] that moreover there holds
\begin{equation}
 K(\Oo_G) \subset K(N)
\end{equation}
for all open neighbourhoods $N$ (in $M$) of $G$, whenever 
$\Oo_G$ is a diamond. Using this, one obtains that the map
$(K,\kappa) \to (S,\sigma)$ given by $[f] \mapsto Ef$ is
surjective, and by Lemma A.1 in [16], it is even a
symplectomorphism. Clearly,
 $(K(\Oo_G),\kappa|K(\Oo_G))$ is a symplectic subspace
of $(K,\kappa)$ for each diamond $\Oo_G$
in $(M,g)$.  For any such diamond one
then obtains, upon viewing it
(or its connected components separately), equipped with the appropriate
restriction of the spacetime metric $g$, as a globally
hyperbolic spacetime in its own right, local versions
of the just introduced symplectic spaces and the symplectomorphisms
between them. More precisely, if we denote by $S(\Oo_G)$
the set of all smooth solutions of the KG-equation (3.3)
with the property that their Cauchy-data on $\Sigma$ are
compactly supported in $G$, then the map
$P_{\Sigma}$ restricts to a symplectomorphism $(S(\Oo_G),
\sigma|S(\Oo_G)) \to ({\cal D}_{G},\delta_{G})$,
$\varphi \mapsto P_{\Sigma}(\varphi)$. Likewise, the
symplectomorphism $[f] \mapsto Ef$ restricts to a symplectomorphism
$(K(\Oo_G),\kappa|K(\Oo_G)) \to
(S(\Oo_G),\sigma|S(\Oo_G))$. 

To the symplectic space $(K,\kappa)$ we can now associate its
Weyl-algebra $\AA[K,\kappa]$, cf.\ Chapter 2. Using the
afforementioned local net-structure of the symplectic space
$(K,\kappa)$, one arrives at the following result.
\begin{Proposition} 
{\rm [16]}. Let $(M,g)$ be a globally hyperbolic
spacetime, and $(K,\kappa)$ the symplectic space, constructed
as above, for the KG-eqn.\ with smooth potential term $r$ on $(M,g)$.
Its Weyl-algebra $\AA[K,\kappa]$ will be called the {\em Weyl-algebra
of the KG-field with potential term $r$ over} $(M,g)$. Define
for each open, relatively compact $\Oo \subset M$, the set $\AA(\Oo)$
as the $C^*$-subalgebra of $\AA[K,\kappa]$ generated by all the
Weyl-operators $W([f])$, $[f] \in K(\Oo)$. Then
$\Oo \to \AA(\Oo)$ is a net of $C^*$-algebras fulfilling
isotony (3.1) and locality (3.2), and moreover {\em primitive
causality}, i.e.\ 
\begin{equation}
\AA(\Oo_G) \subset \AA(N)
\end{equation}
for all neighbourhoods $N$ (in $M$) of $G$, whenever $\Oo_G$ is
a (relatively compact) diamond.
\end{Proposition}
It is worth recalling (cf.\ [5]) that the Weyl-algebras
corresponding to symplectically equivalent spaces are
canonically isomorphic in the following way: Let
$W(x)$, $x \in K$ denote the Weyl-generators of $\AA[K,\kappa]$
and $W_S(\varphi)$, $\varphi \in S$, the Weyl-generators
of $\AA[S,\sigma]$. Furthermore, let $T$ be a symplectomorphism
between $(K,\kappa)$ and $(S,\sigma)$.
                    Then there is a uniquely determined
$C^*$-algebraic isomorphism $\alpha_T : \AA[K,\kappa] \to
\AA[S,\sigma]$ given by $\alpha_T(W(x)) = W_S(Tx)$, $x \in K$.
This shows that if we had associated e.g.\ with $(S,\sigma)$
the Weyl-algebra $\AA[S,\sigma]$ as the algebra of quantum observables
of the KG-field over $(M,g)$, we would have obtained an
equivalent net of observable algebras
(connected to the previous one by a net isomorphism,
see [3,16]), rendering the same
physical information.
\\[24pt]
{\bf 3.4 Hadamard States}
\\[18pt]
We have indicated above that quasifree Hadamard states
are distinguished by their short-distance behaviour which
allows the definition of expectation values of energy-momentum
observables with reasonable properties [26,68,69,71]. If
$\omega_{\mu}$ is a quasifree state on the Weyl-algebra
$\AA[K,\kappa]$, then we call
$$ \lambda(x,y) := \mu(x,y) + \frac{i}{2}\kappa(x,y)\,, \quad
x,y \in K\,, $$
its {\it two-point function} and
$$ \Lambda(f,h) := \lambda([f],[h])\,, \quad f,h \in \CoM\,, $$
its {\it spatio-temporal} two-point function.
In Chapter 2 we have seen that a quasifree state is entirely
determined through specifying $\mu \in {\sf q}(K,\kappa)$, which
is equivalent to the specification of the two-point function
$\lambda$. Sometimes the notation $\lambda_{\omega}$
or $\lambda_{\mu}$ will be used to indicate the quasifree
state $\omega$ or the dominating scalar product $\mu$ which is
determined by $\lambda$.

For a quasifree Hadamard state, the spatio-temporal two-point
function is of a special form, called Hadamard form. The
definition of Hadamard form which we give here follows that due to
Kay and Wald [45].  Let $N$ is a causal normal neighbourhood
of a Cauchy-surface $\Sigma$ in $(M,g)$. Then a smooth function
$\chi : N \times N \to [0,1]$ is called {\it $N$-regularizing}
if it has the following property: There is an open
neighbourhood, $\Omega_*$, in $N \times N$ of the set of
pairs of causally related points in $N$ such that
$\overline{\Omega_*}$ is contained in a set $\Omega$ to
be described presently, and $\chi \equiv 1$ on $\Omega_*$
while $\chi \equiv 0$ outside of $\overline{\Omega}$. Here,
$\Omega$ is an open neighbourhood in $M \times M$ of the
set of those $(p,q) \in M \times M$ which are causally related
and have the property that (1) $J^+(p) \cap J^-(q)$ and
$J^+(q) \cap J^-(p)$ are contained within a convex normal
neighbourhood, and (2) $s(p,q)$, the square of the geodesic
distance between $p$ and $q$, is a well-defined, smooth
function on $\Omega$. (One observes that there are always
sets $\Omega$ of this type which contain a neighbourhood of the
diagonal in $M \times M$, and that an $N$-regularizing function
depends on the choice of the pair of sets $\Omega_*,\Omega$
with the stated properties.) It is not difficult to check that
$N$-regularizing functions always exist for any causal normal
neighbourhood; a proof of that is e.g.\ given in [55].
Then denote by $U$ the square root of the VanVleck-Morette determinant,
and by $v_m$, $m \in \NN_0$ the sequence determined by the
Hadamard recursion relations for the KG-equation (3.3),
see [23,27] and also [30] for their definition.
They are all smooth functions on $\Omega$.\footnote{For any
choice of $\Omega$ with the properties just described.}
Now set for $n \in \NN$,
$$ V^{(n)}(p,q) := \sum_{m = 0}^n v_m(p,q)(s(p,q))^m \,,
\quad (p,q) \in \Omega\,, $$
and, given a smooth time-function $T: M \to \RR$ increasing
towards the future, define for all $\epsilon > 0$ and
$(p,q) \in \Omega$,
$$ Q_T(p,q;\epsilon) := s(p,q) - 2 i\epsilon (T(p) - T(q)) -
\epsilon^2 \,,$$
and
$$G^{T,n}_{\epsilon}(p,q) := \frac{1}{4\pi^2}\left(
 \frac{U(p,q)}{Q_T(p,q;\epsilon)} + V^{(n)}(p,q)ln(Q_T(p,q;
 \epsilon)) \right) \,, $$
 where $ln$ is the principal branch of the logarithm.
With this notation, one can give the
\begin{Definition}{\rm [45]}. A $\CC$-valued
bilinear form $\Lambda$ on $\CoM$ is called an {\em Hadamard form}
if, for a suitable choice of a causal normal neighbourhood
$N$ of some Cauchy-surface $\Sigma$, and for  suitable choices  of an
$N$-regularizing function $\chi$ and a future-increasing
time-function $T$ on $M$, there exists a sequence
$H^{(n)} \in C^n(N \times N)$, so that
\begin{equation}
\Lambda(f,h) = \lim_{\epsilon \to 0+}
\int_{M\times M} \Lambda^{T,n}_{\epsilon}(p,q)f(p)h(q)\,
d\eta(p) \,d\eta(q)
\end{equation}
for all $f,h \in C_0^{\infty}(N,\RR)$, where
\footnote{ The set $\Omega$ on which the functions forming
$G^{T,n}_{\epsilon}$ are defined and smooth is here to
coincide with the $\Omega$ with respect to which $\chi$ is
defined.}
\begin{equation}
\Lambda^{T,n}_{\epsilon}(p,q) := \chi(p,q)G^{T,n}_{\epsilon}(p,q)
 + H^{(n)}(p,q)\,,
\end{equation}
and if, moreover, $\Lambda$ is a global bi-parametrix of
the KG-equation (3.3), i.e.\ it satisfies
$$ \Lambda((\nabla^a\nabla_a + r)f,h) = B_1(f,h)\quad {\it and}
\quad \Lambda(f,(\nabla^a\nabla_a + r)h) = B_2(f,h) $$
for all $f,h \in C_0^{\infty}(M)$, where $B_1$ and $B_2$ are
given by smooth integral kernels on $M \times M$.\footnote{
We point out that statement (b) of Prop.\ 3.4 is wrong if
the assumption that $\Lambda$ is a global bi-parametrix is not made.
In this respect, Def.\ C.1 of [66] is imprecisely formulated as
the said assumption is not stated. There, like in several other
references, it has been implicitely assumed that $\Lambda$ is a
two point function and thus a bi-solution of (3.3), i.e. a
bi-parametrix with $B_1 = B_2 \equiv 0$.}
\end{Definition}
Based on results of [24,25], it is shown in [45] that this is
a reasonable definition. The findings of these works will be
collected in the following
\begin{Proposition} ${}$\\[6pt]
(a) If $\Lambda$ is of Hadamard form on a causal normal neighbourhood
$ N$ of a Cauchy-surface $\Sigma$ for some choice of a time-function
$T$ and some  $N$-regularizing function $\chi$ (i.e.\ that
(3.9),(3.10) hold with suitable $H^{(n)} \in C^n(N \times
N)$), then so it is for any other time-function $T'$ and
$N$-regularizing $\chi'$. (This means that these changes
can be compensated by choosing another sequence $H'^{(n)} \in
C^n( N \times N)$.)
\\[6pt]
(b) (Causal Propagation Property of the Hadamard Form)\\
If $\Lambda$ is of Hadamard form on a causal normal neighbourhood 
$ N$ of some Cauchy-surface $\Sigma$, then it is of Hadamard form
in any causal normal neighbourhood $ N'$ of any other Cauchy-surface
$\Sigma'$.
\\[6pt]
(c) Any $\Lambda$ of Hadamard form is a regular kernel distribution
on $C_0^{\infty}(M \times M)$.
\\[6pt]
(d) There exist pure, quasifree Hadamard states (these will be
referred to as {\em Hadamard vacua}) on the Weyl-algebra $\AA[K,\kappa]$ of the 
KG-field in any globally hyperbolic spacetime. The family of quasifree 
Hadamard states on $\AA[K,\kappa]$ spans an infinite-dimensional subspace 
of the continuous dual space of $\AA[K,\kappa]$.
\\[6pt]
(e) The dominating scalar products $\mu$ on $K$ arising from quasifree 
Hadamard states $\omega_{\mu}$ induce locally the same topology,
i.e.\ if $\mu$ and $\mu'$ are arbitrary such scalar products and
$\Oo \subset M$ is open and relatively compact, then there are two
positive constants $a,a'$ such that
$$ a\, \mu([f],[f]) \leq \mu'([f],[f]) \leq a'\,\mu([f],[f])\,,
\quad [f] \in K(\Oo)\,.$$
\end{Proposition}
{\it Remark.} Observe that this definition of Hadamard form rules out
the occurence of spacelike singularities, meaning that the Hadamard form
$\Lambda$ is, when tested on functions $f,h$ in (3.9) whose
supports are acausally separated, given by a $C^{\infty}$-kernel.
For that reason, the definition of Hadamard form as stated
above is also called {\it global} Hadamard form (cf.\ [45]).
A weaker definition of Hadamard form would be to prescribe (3.9),(3.10)
only for sets $N$ which, e.g., are members of an open covering of $M$
by convex normal neighbourhoods, and thereby to require the Hadamard form
locally. In the case that $\Lambda$ is the spatio-temporal two-point
function of a state on $\AA[K,\kappa]$ and thus dominates the symplectic
form $\kappa$ ($|\kappa([f],[h])|^2 \leq 4\,\Lambda(f,f)\Lambda(h,h)$),
it was recently proved by Radzikowski that if $\Lambda$ is locally
of Hadamard form, then it is already globally of Hadamard form [56].
However, if $\Lambda$ doesn't dominate $\kappa$, this need not hold
[29,51,56]. Radzikowski's proof makes use of a characterization
of Hadamard forms in terms of their wave front sets which was mentioned
above. A definition of Hadamard form which is less technical in appearence
has recently been given in [44].
We should add that the usual Minkowski-vacuum of the free scalar
field with constant, non-negative potential
 term is, of course, an Hadamard vacuum.
This holds, more generally, also for ultrastatic spacetimes, see below.
\\[10pt]
{\it Notes on the proof of Proposition 3.4.}
The property (a) is proved in [45]. The argument for (b) is 
essentially contained in [25] and in the generality stated here it is
completed in [45]. An alternative proof using the 
``propagation of singularities theorem'' for hyperbolic differential 
equations is presented in [55].
 Also property (c) is proved in [45 (Appendix B)] (cf.\
[66 (Prop.\ C.2)]). The existence of Hadamard vacua (d)
is proved in [24] (cf.\ also [45]);  the stated Corollary has
been observed in [66] (and, in slightly different formulation,
already in [24]). Statement (e) has been shown to hold in
[66 (Prop.\ 3.8)].
\\[10pt]
In order to prepare the formulation of the next result, in which we will
apply our result of Chapter 2, we need to collect some more notation.
 Suppose that we are given a quasifree state $\omega_{\mu}$ on
the Weyl-algebra $\AA[K,\kappa]$ of the KG-field over some
globally hyperbolic spacetime $(M,g)$, and that $\Sigma$ is a
Cauchy-surface in that spacetime. Then we denote by $\mu_{\Sigma}$
the dominating scalar product on $({\cal D}_{\Sigma},\delta_{\Sigma})$
which is, using the symplectomorphism between $(K,\kappa)$ and
$({\cal D}_{\Sigma},\delta_{\Sigma})$, induced by the dominating scalar
product $\mu$ on $(K,\kappa)$, i.e.\
\begin{equation}
 \mu_{\Sigma}(P_{\Sigma}Ef,P_{\Sigma}Eh) = \mu([f],[h])\,, \quad
[f],[h] \in K\,. 
\end{equation}
Conversely, to any $\mu_{\Sigma} \in {\sf q}(\DS,\delta_{\Sigma})$
there corresponds via (3.11) a $\mu \in {\sf q}(K,\kappa)$.

Next, consider a complete Riemannian manifold $(\Sigma,\gamma)$, with 
corresponding Laplacian $\Delta_{\gamma}$, and as before, consider the
operator $ -\Delta_{\gamma} +1$ on $C_0^{\infty}(\Sigma,\RR)$.
 Owing to the completeness
of $(\Sigma,\gamma)$ this operator is,
together with all its powers, essentially selfadjoint in
$L^2_{\RR}(\Sigma,d\eta_{\gamma})$ [10],
and we denote its selfadjoint extension
by $A_{\gamma}$. Then one can introduce the
{Sobolev scalar products} of $m$-th order,
$$ \langle u,v \rangle_{\gamma,m} := \langle u, A_{\gamma}^m v \rangle\,,
\quad u,v \in C_0^{\infty}(\Sigma,\RR),\  m \in \RR\,, $$
where on the right hand side is the scalar product of $L^2_{\RR}(\Sigma,
d\eta_{\gamma})$. The completion of $C_0^{\infty}(\Sigma,\RR)$
in the topology of $\langle\,.\,,\,.\,\rangle_{\gamma,m}$
will be denoted by $H_m(\Sigma,\gamma)$.
It turns out that the topology of $H_m(\Sigma,\gamma)$ is locally
independent of the complete Riemannian metric $\gamma$, and that
composition with diffeomorphisms and multiplication with smooth,
compactly supported functions are continuous operations on these
Sobolev spaces. (See Appendix B for precise formulations of these
statements.) Therefore, whenever $G \subset \Sigma$ is open and
relatively compact, the topology which $\langle \,.\,,\,.\, \rangle_{m,\gamma}$
induces on $C_0^{\infty}(G,\RR)$ is independent of the particular
complete Riemannian metric $\gamma$, and we shall refer to the
topology which is thus locally induced on $C_0^{\infty}(\Sigma,\RR)$
simply as the (local) {\it $H_m$-topology.}

Let us now suppose that we have an ultrastatic spacetime $(\MT,\gT)$,
given in a natural foliation as $(\RR \times \ST,dt^2 \otimes (-\gamma))$
where $(\ST,\gamma)$ is a complete Riemannian manifold. We shall
identify $\ST$ and $\{0\} \times \ST$. Consider again
$A_{\gamma}$  = selfadjoint extension of $- \Delta_{\gamma} + 1$ on
$C_0^{\infty}(\ST,\RR)$ in $L^2_{\RR}(\ST,d\eta_{\gamma})$ with
$\Delta_{\gamma}$ = Laplacian of $(\ST,\gamma)$, and the scalar product
$\mu^{\circ}_{\ST}$ on ${\cal D}_{\ST}$ given by
\begin{eqnarray}
 \mu^{\circ}_{\ST}(\uu,\vv) & := & \frac{1}{2} \left ( 
                       \langle u_0,A_{\gamma}^{1/2}v_0 \rangle
                  + \langle u_1,A_{\gamma}^{-1/2}v_1 \rangle \right) \\
 & = & \frac{1}{2} \left(  \langle u_0,v_0 \rangle_{\gamma,1/2}
 + \langle u_1,v_1 \rangle_{\gamma,-1/2} \right) \nonumber
\end{eqnarray}
for all $\uu,\vv \in {\cal D}_{\ST}$. It is now straightforward to
check that $\mu^{\circ}_{\ST} \in {\sf pu}({\cal D}_{\ST},\delta_{\ST})$,
 in fact, $\mu^{\circ}_{\ST}$ is the purification of the classical energy
scalar product $\mu^E_{\ST}$ defined in eqn.\ (3.5). (We refer to
[11] for discussion, and also the treatment of more general situations
along similar lines.) What is furthermore central for the derivation of
the next result is that $\mu^{\circ}_{\ST}$ corresponds (via (3.11)) to
 an Hadamard vacuum $\omega^{\circ}$ on the Weyl-algebra
of the KG-field with potential term $r \equiv 1$ over the ultrastatic
spacetime $(\RR \times \ST,dt^2 \oplus (-\gamma))$. This has been
proved in [24]. The state $\omega^{\circ}$ is called the
{\it ultrastatic vacuum} for the said KG-field over 
$(\RR \times \ST ,dt^2 \oplus (-\gamma))$; it is the unique pure, 
quasifree ground state on the corresponding Weyl-algebra for the
time-translations $(t,q) \mapsto (t + t',q)$ on that ultrastatic
spacetime with respect to the chosen natural foliation (cf.\ [40,42]).
\\[6pt]
{\it Remark.} The passage from $\mu^E_{\ST}$ to $\mu^{\circ}_{\ST}$,
where $\mu^{\circ}_{\ST}$ is the purification of the classical
energy scalar product, may be viewed as a refined form of
``frequency-splitting'' procedures (or Hamiltonian diagonalization),
in order to obtain pure dominating scalar products and hence, pure states
of the KG-field in curved spacetimes, see [11]. However, in the case
that $\ST$ is not a Cauchy-surface lying in the natural foliation of
an ultrastatic spacetime, but an arbitrary Cauchy-surface in an
arbitrary globally hyperbolic spacetime, the $\mu^{\circ}_{\ST}$
may fail to correspond to a quasifree Hadamard state --- even though,
as the following Proposition demonstrates, $\mu^{\circ}_{\ST}$ gives
locally on the Cauchy-data space ${\cal D}_{\ST}$ the same topology
as the dominating scalar products induced on it by any quasifree
Hadamard state. More seriously, $\mu^{\circ}_{\ST}$ may even correspond to
a state which is no longer locally quasiequivalent to any
quasifree Hadamard state. For an explicit example demonstrating
this in a closed Robertson-Walker universe, and for additional
discussion, we refer to Sec.\ 3.6 in [38].
\\[6pt] 
We shall say that a map $T : {\cal D}_{\Sigma} \to {\cal D}_{\Sigma'}$,
with $\Sigma,\Sigma'$ Cauchy-surfaces, is {\it locally continuous} if,
for any open, locally compact $G \subset \Sigma$, the restriction of
$T$ to $C_0^{\infty}(G,\RR) \oplus C_0^{\infty}(G,\RR)$ is continuous
(with respect to the topologies under consideration).
\begin{Proposition}
Let $\omega_{\mu}$ be a quasifree Hadamard state on the Weyl-algebra
$\AA[K,\kappa]$ of the KG-field with smooth potential term $r$ over the
globally hyperbolic spacetime $(M,g)$, and $\Sigma,\Sigma'$ two 
Cauchy-surfaces in $(M,g)$.

Then the Cauchy-data evolution map
\begin{equation}
 T_{\Sigma',\Sigma} : = P_{\Sigma'} \lcrc P_{\Sigma}^{-1} :
 {\cal D}_{\Sigma} \to {\cal D}_{\Sigma'}
\end{equation}
is locally continuous in the $H_{\tau} \oplus H_{\tau -1}$-topology,
$0 \leq \tau \leq 1$, on the Cauchy-data spaces, and the topology
induced by $\mu_{\Sigma}$ on ${\cal D}_{\Sigma}$ coincides locally
(i.e.\ on each $C_0^{\infty}(G,\RR) \oplus C_0^{\infty}(G,\RR)$
for $G \subset \Sigma$ open and relatively compact) with the
 $H_{1/2} \oplus H_{-1/2}$-topology.
\end{Proposition}
{\it Remarks.} (i) Observe that the continuity statement is
reasonably formulated since, as a consequence of the support
properties of solutions of the KG-equation with Cauchy-data
of compact support (``finite propagation speed'') it holds that
for each open, relatively compact $G \subset \Sigma$ there is 
an open, relatively compact $G' \subset \Sigma'$ with
$T_{\Sigma',\Sigma}(C_0^{\infty}(G,\RR) \oplus C_0^{\infty}(G,\RR))
\subset C_0^{\infty}(G',\RR) \oplus C_0^{\infty}(G',\RR)$.
 \\[6pt]
(ii) For $\tau =1$, the continuity statement is just the
classical energy estimate.
 It should be mentioned here that the claimed continuity can
also be obtained by other methods. For instance, Moreno [50]
proves, under more restrictive assumptions on $\Sigma$ and $\Sigma'$
(among which is their compactness), the continuity of $T_{\Sigma',\Sigma}$
in the topology of $H_{\tau} \oplus H_{\tau -1}$ for all $\tau \in \RR$,
by employing an abstract energy estimate for first order hyperbolic equations
(under suitable circumstances, the KG-equation can be brought into this form).
We feel, however, that our method, using the results of Chapter 2, is
physically more appealing and emphasizes much better the ``invariant''
structures involved, quite in keeping with the general approach to quantum
field theory. 
\\[10pt]
{\it Proof of Proposition 3.5.} We note that there is a diffeomorphism
$\Psi : \Sigma \to \Sigma'$. To see this, observe that we may pick
a foliation $F : \RR \times \ST \to M$ of $M$ in Cauchy-surfaces. Then
for each $q \in \ST$, the curves $t \mapsto F(t,q)$ are inextendible,
timelike curves in $(M,g)$. Each such curve intersects $\Sigma$ exactly
once, at the parameter value $t = \tau(q)$. Hence $\Sigma$ is the set
$\{F(\tau(q),q) : q \in \ST\}$. As $F$ is a diffeomorphism and
$\tau: \ST \to \RR$ must be $C^{\infty}$ since, by assumption,
$\Sigma$ is a smooth hypersurface in $M$, one can see that $\Sigma$
and $\ST$ are diffeomorphic. The same argument shows that 
$\Sigma'$ and $\ST$ and therefore, $\Sigma$ and $\Sigma'$, are
diffeomorphic.

Now let us first assume that the $g$-induced Riemannian metrics
$\gamma_{\Sigma}$ and $\gamma_{\Sigma'}$ on $\Sigma$, resp.\
$\Sigma'$, are complete. Let $d\eta$ and $d\eta'$ be the induced
volume measures on $\Sigma$ and $\Sigma'$, respectively. The $\Psi$-transformed
measure of $d\eta$ on $\Sigma'$, $\Psi^*d\eta$, is given through
\begin{equation}
 \int_{\Sigma} (u \lcrc \Psi) \,d\eta = \int_{\Sigma'} u\,(\Psi^*d\eta)\,,
\quad u \in C_0^{\infty}(\Sigma')\,.
\end{equation}
Then the Radon-Nikodym derivative $(\rho(q))^2 :=(\Psi^*d\eta/d\eta')(q)$,
$q \in \Sigma'$, is a smooth, strictly positive function on $\Sigma'$,
and it is now easy to check that the linear map
$$ \vartheta : ({\cal D}_{\Sigma},\delta_{\Sigma})
\to ({\cal D}_{\Sigma'},\delta_{\Sigma'})\,, \quad
\uu \mapsto \rho \cdot (u_0 \lcrc \Psi^{-1}) \oplus \rho \cdot
 (u_1 \lcrc \Psi^{-1}) \,, $$
is a symplectomorphism. Moreover, by the result given in Appendix
B, $\vartheta$ and its inverse are locally continuous maps in the
$H_s \oplus H_t$-topologies on both Cauchy-data spaces, for all
$s,t \in \RR$.

By the energy estimate, $T_{\Sigma',\Sigma}$ is locally continuous 
with respect to the $H_1 \oplus H_0$-topology on the Cauchy-data
spaces, and the same holds for the inverse $(T_{\Sigma',\Sigma})^{-1}
= T_{\Sigma,\Sigma'}$. Hence, the map 
$\Theta := \vartheta^{-1} \lcrc T_{\Sigma',\Sigma}$ is a symplectomorphism
of $(\DS,\delta_{\Sigma})$, and $\Theta$ together with its inverse
is locally continuous in the $H_1 \oplus H_0$-topology on $\DS$.
Here we made use of Remark (i) above. Now pick two sets $G$ and $G'$
as in Remark (i), then there is some open, relatively compact neighbourhood
$\tilde{G}$ of $\Psi^{-1}(G') \cup G$ in $\Sigma$. We can choose a smooth,
real-valued function $\chi$ compactly supported on $\Sigma$ with $\chi \equiv
 1$ on $\tilde{G}$. It is then straightforward to check that the maps
$\chi \lcrc \Theta \lcrc \chi$ and $\chi \lcrc \Theta^{-1} \lcrc \chi$
($\chi$ to be interpreted as multiplication with $\chi$) is a pair
of symplectically adjoint maps on $(\DS,\delta_{\Sigma})$ which are bounded
with respect to the $H_1 \oplus H_0$-topology, i.e.\ with respect to the
norm of $\mu_{\Sigma}^E$. At this point we use Theorem 2.2(b) and consequently
$\chi \lcrc \Theta \lcrc \chi$ and $\chi \lcrc \Theta^{-1}\lcrc \chi$ are
continuous with respect to the norms of the $(\mu^E_{\Sigma})_s$,
$0 \leq s \leq 2$. Inspection shows that
$$ (\mu^E_{\Sigma})_s (\uu,\vv) =
\frac{1}{2} \left( \langle u_0,A_{\gamma_{\Sigma}}^{1-s/2}v_0 \rangle
 + \langle u_1,A_{\gamma_{\Sigma}}^{-s/2}v_1 \rangle \right)
$$
for $0 \leq s \leq 2$. From this it is now easy to see that
$\Theta$ restricted to $C_0^{\infty}(G,\RR) \oplus C_0^{\infty}(G,\RR)$
is continuous in the topology of $H_{\tau} \oplus H_{\tau -1}$,
$0 \leq \tau \leq 1$, since 
$\chi \lcrc \Theta \lcrc \chi(\uu) = \Theta(\uu)$  for all
$\uu \in C_0^{\infty}(G,\RR) \oplus C_0^{\infty}(G,\RR)$
by the choice of $\chi$. Using that $\Theta = \vartheta^{-1}\lcrc T_{
\Sigma',\Sigma}$ and that $\vartheta$ is locally continuous with respect to 
all the $H_s \oplus H_t$-topologies, $s,t \in \RR$, on the Cauchy-data
spaces, we deduce that that $T_{\Sigma',\Sigma}$ is locally continuous
in the $H_{\tau} \oplus H_{\tau -1}$-topology, $0 \leq \tau \leq 1$,
as claimed.

If the $g$-induced Riemannian metrics $\gamma_{\Sigma}$, $\gamma_{\Sigma'}$
are not complete, one can make them into complete ones
$\hat{\gamma}_{\Sigma} := f \cdot \gamma_{\Sigma}$, $\hat{\gamma}_{\Sigma'}
:= h \cdot \gamma_{\Sigma'}$ by multiplying them with suitable smooth,
strictly positive functions $f$ on $\Sigma$ and $h$ on $\Sigma'$ [14].
Let $d\hat{\eta}$ and $d\hat{\eta}'$ be the volume measures corresponding 
to the new metrics. Consider then the density functions 
$(\phi_1)^2 := (d\eta/d\hat{\eta})$,
 $(\phi_2)^2 := (d\hat{\eta}'/d\eta')$,
which are $C^{\infty}$ and strictly positive, and define 
$(\DS,\hat{\delta}_{\Sigma})$, $({\cal D}_{\Sigma'},\hat{\delta}_{\Sigma'})$
and $\hat{\vartheta}$ like their unhatted counterparts but with
$d\hat{\eta}$ and $d\hat{\eta}'$ in place of $d\eta$ and $d\eta'$.
Likewise define $\hat{\mu}^E_{\Sigma}$ with respect to $\hat{\gamma}_{\Sigma}$.
Then $\hat{T}_{\Sigma',\Sigma} := \phi_2 \lcrc T_{\Sigma',\Sigma}
\lcrc \phi_1$ (understanding that $\phi_1,\phi_2$ act as
multiplication operators) and its inverse are symplectomorphisms
between $(\DS,\hat{\delta}_{\Sigma})$ and
$({\cal D}_{\Sigma'},\hat{\delta}_{\Sigma'})$ which are  locally
continuous in the $H_1 \oplus H_0$-topology. Now we can apply the
argument above showing that $\hat{\Theta} = \hat{\vartheta}^{-1} \lcrc
\hat{T}_{\Sigma',\Sigma}$ and, hence, $\hat{T}_{\Sigma',\Sigma}$ is
locally continuous in the $H_{\tau} \oplus H_{\tau -1}$-topology for
$0 \leq \tau \leq 1$. The same follows then for 
$T_{\Sigma',\Sigma} = \phi_2^{-1} \lcrc \hat{T}_{\Sigma',\Sigma}
\lcrc \phi_1^{-1}$.

For the proof of the second part of the statement, we note first
that in [24] it is shown that there exists another globally
hyperbolic spacetime $(\hat{M},\hat{g})$ of the form
$\hat{M} = \RR \times \Sigma$ with the following properties:
\\[6pt]
(1) $\Sigma_0 : = \{0\} \times \Sigma$ is a Cauchy-surface in $(\hat{M},
\hat{g})$, and a causal normal neighbourhood $N$ of $\Sigma$ in $M$
coincides with a causal normal neigbourhood $\hat{N}$ of
$\Sigma_{0}$ in $\hat{M}$, in such a way that $\Sigma = \Sigma_0$
and $g = \hat{g}$ on $N$.
\\[6pt]
(2) For some $t_0 < 0$, the $(-\infty,t_0) \times \Sigma$-part of
$\hat{M}$ lies properly to the past of $\hat{N}$, and on that part,
$\hat{g}$ takes the form $dt^2 \oplus (- \gamma)$ where
$\gamma$ is a complete Riemannian metric on $\Sigma$.
\\[6pt]
This means that $(\hat{M},\hat{g})$ is a globally hyperbolic
spacetime which equals $(M,g)$ on a causal normal neighbourhood of $\Sigma$
and becomes ultrastatic to the past of it.

Then consider the Weyl-algebra $\AA[\hat{K},\hat{\kappa}]$ of the 
KG-field with potential term $\hat{r}$ over $(\hat{M},\hat{g})$, where 
$\hat{r} \in C_0^{\infty}(\hat{M},\RR)$ agrees with $r$ on the
neighbourhood $\hat{N} = N$ and is identically equal to $1$ on the
$(-\infty,t_0) \times \Sigma$-part of $\hat{M}$. Now observe that
the propagators $E$ and $\hat{E}$ of the respective KG-equations
on $(M,g)$ and $(\hat{M},\hat{g})$ coincide when restricted to
$C_0^{\infty}(N,\RR)$. Therefore one obtains an identification map 
$$ [f] = f + {\rm ker}(E) \mapsto [f]\,\hat{{}} = f + {\rm ker}(\hat{E}) \,,
\quad f \in C_0^{\infty}(N,\RR) \,,$$
between $K(N)$ and $\hat{K}(\hat{N})$ which preserves the 
symplectic forms $\kappa$ and $\hat{\kappa}$. Without danger we may
write this identification as an equality,
$K(N) = \hat{K}(\hat{N})$.
This identification map between $(K(N),\kappa|K(N))$
and $(\hat{K}(\hat{N}),\hat{\kappa}|\hat{K}(\hat{N}))$ lifts to
a $C^*$-algebraic isomorphism between the corresponding
Weyl-algebras
\begin{eqnarray}
 \AA[K(N),\kappa|K(N)]& =& \AA[\hat{K}(\hat{N}),\hat{\kappa}|
\hat{K}(\hat{N})]\,, \nonumber \\
 W([f])& =& \hat{W}([f]\,\hat{{}}\,)\,,\ \ \
f \in C_0^{\infty}(N,\RR)\,.
\end{eqnarray}
Here we followed our just indicated convention to abbreviate
this identification as an equality. Now we have
$D(N) = M$ in $(M,g)$ and $D(\hat{N}) = \hat{M}$ in
$(\hat{M},\hat{g})$, implying that $K(N) = K$ and
$\hat{K}(\hat{N}) = \hat{K}$. Hence $\AA[K(N),\kappa|K(N)] =
\AA[K,\kappa]$ and the same for the ``hatted'' objects.
Thus (3.15) gives rise to an identification between
$\AA[K,\kappa]$ and $\AA[\hat{K},\hat{\kappa}]$, and so the
quasifree Hadamard state $\omega_{\mu}$ induces a quasifree
state $\omega_{\hat{\mu}}$ on $\AA[\hat{K},\hat{\kappa}]$
with
\begin{equation}
\hat{\mu}([f]\,\hat{{}},[h]\,\hat{{}}\,) = \mu([f],[h])\,, \quad
f,h \in C_0^{\infty}(N,\RR) \,.
\end{equation}
This state is also an Hadamard state since we have
\begin{eqnarray*}
\Lambda(f,h)& =& \mu([f],[h]) + \frac{i}{2}\kappa([f],[h]) \\
 & = & \hat{\mu}([f]\,\hat{{}}\,,[h]\,\hat{{}}\,) + \frac{i}{2}
\hat{\kappa}([f]\,\hat{{}}\,,[h]\,\hat{{}}\,)\,, \quad f,h \in
C_0^{\infty}(N,\RR)\,,
\end{eqnarray*}
and $\Lambda$ is, by assumption, of Hadamard form.
However, due to the causal propagation property of the
Hadamard form this means that $\hat{\mu}$ is the dominating
scalar product on $(\hat{K},\hat{\kappa})$ of a quasifree
Hadamard state on $\AA[\hat{K},\hat{\kappa}]$.
Now choose some $t < t_0$, and let $\Sigma_t = \{t\} \times
\Sigma$ be the Cauchy-surface in the ultrastatic part
of $(\hat{M},\hat{g})$ corresponding to this value of the
time-parameter of the natural foliation. As remarked above,
the  scalar product
\begin{equation}
 \mu^{\circ}_{\Sigma_t}(\uu,\vv)
 = \frac{1}{2}\left( \langle u_0,v_0 \rangle_{\gamma,1/2}
 + \langle u_1,v_1 \rangle_{\gamma,-1/2} \right)\,,
\quad \uu,\vv \in {\cal D}_{\Sigma_t} \,,
\end{equation}
is the dominating scalar product on
$({\cal D}_{\Sigma_t},\delta_{\Sigma_t})$ corresponding to the
ultrastatic vacuum state $\omega^{\circ}$ over the 
ultrastatic part of  $(\hat{M},\hat{g})$, which is an Hadamard
vacuum. Since the dominating scalar products of all
quasifree Hadamard states yield locally the same topology
(Prop.\ 3.4(e)), it follows that the dominating scalar product
$\hat{\mu}_{\Sigma_t}$ on $({\cal D}_{\Sigma_{t}},\delta_{\Sigma_t})$,
which is induced (cf.\ (3.11)) by the the dominating scalar product
of $\hat{\mu}$ of the quasifree Hadamard state $\omega_{\hat{\mu}}$,
endows ${\cal D}_{\Sigma_t}$ locally with the same topology
as does $\mu^{\circ}_{\Sigma_t}$. As can be read off from (3.17),
this is the local $H_{1/2} \oplus H_{-1/2}$-topology.

To complete the argument, we note that (cf.\ (3.11,3.13))
$$ \hat{\mu}_{\Sigma_0}(\uu,\vv) =
\hat{\mu}_{\Sigma_t}(T_{\Sigma_t,\Sigma_0}(\uu),T_{\Sigma_t,\Sigma_0}
(\vv))\,, \quad \uu,\vv \in {\cal D}_{\Sigma_0}\,.$$
But since $\hat{\mu}_{\Sigma_t}$ induces locally the
$H_{1/2} \oplus H_{-1/2}$-topology and since the symplectomorphism
$T_{\Sigma_t,\Sigma_0}$ as well as its inverse are locally continuous
on the Cauchy-data spaces in the $H_{1/2}\oplus H_{-1/2}$-topology,
the last equality entails that $\hat{\mu}_{\Sigma_0}$ induces the
local $H_{1/2} \oplus H_{-1/2}$-topology on ${\cal D}_{\Sigma_0}$.
In view of (3.16), the Proposition is now proved. $\Box$   
\\[24pt]
{\bf 3.5 Local Definiteness, Local Primarity,
Haag-Duality, etc.}
\\[18pt]
In this section we prove Theorem 3.6 below on the algebraic
structure of the GNS-representations associated with quasifree
Hadamard states on the CCR-algebra of the KG-field on an
arbitrary globally hyperbolic spacetime $(M,g)$. The
results appearing therein extend our previous work [64,65,66].

Let $(M,g)$ be a globally hyperbolic spacetime. 
We recall that a  subset $\Oo$ of $M$
is called a { regular diamond} if it is of the form
$\Oo = \Oo_G = {\rm int}\,D(G)$ where
$G$ is an open, relatively compact subset of some Cauchy-surface
$\Sigma$ in $(M,g)$ having the property that the boundary
$\partial G$ of $G$ is contained in the
union of finitely many smooth, closed, two-dimensional submanifolds 
of $\Sigma$. We also recall the notation ${\cal R}_{\omega}(\Oo)
= \pi_{\omega}(\AA(\Oo))^-$ for the local von Neumann algebras in
the GNS-representation of a state $\omega$. The $C^*$-algebraic
net of observable algebras $\Oo \to \AA(\Oo)$
 will be understood as being that associated
with the KG-field in Prop.\ 3.2.
\begin{Theorem}
Let $(M,g)$ be a globally hyperbolic spacetime and 
$\AA[K,\kappa]$ the Weyl-algebra of the KG-field with smooth, real-valued
potential function $r$ over $(M,g)$. Suppose that 
$\omega$ and $\omega_1$ are two quasifree Hadamard states on
$\AA[K,\kappa]$. Then the following statements hold.
\\[6pt]
(a) The GNS-Hilbertspace $\Hh_{\omega}$
of $\omega$ is infinite dimensional and separable.
\\[6pt]
(b) The restrictions of the GNS-representations $\pi_{\omega}|\AO$
and $\pi_{\omega_1}|\AO$ of any open, relatively compact
$\Oo \subset M$ are quasiequivalent. They are even unitarily
equivalent when $\Oo^{\perp}$ is non-void. 
\\[6pt]
(c) For each $p \in M$ we have local definiteness,
$$ \bigcap_{\Oo \owns p} {\cal R}_{\omega}(\Oo) = \CC \cdot 1\, . $$
More generally, whenever $C \subset M$ is the subset of a compact
set which is contained in the union of finitely many smooth, closed,
two-dimensional submanifolds of an arbitrary Cauchy-surface
$\Sigma$ in $M$,
then 
\begin{equation}
 \bigcap_{\Oo \supset C} {\cal R}_{\omega}(\Oo) = \CC \cdot 1\,.
\end{equation}
\\[6pt]
(d) Let $\Oo$ and $\Oo_1$ be two relatively compact diamonds, based
on Cauchy-surfaces $\Sigma$ and $\Sigma_1$, respectively, such
that $\overline{\Oo} \subset \Oo_1$. Then the split-property
holds for the pair ${\cal R}_{\omega}(\Oo)$ and 
${\cal R}_{\omega}(\Oo_1)$, i.e.\ there exists a type ${\rm I}_{\infty}$
factor $\cal N$ such that one has the inclusion
$$ {\cal R}_{\omega}(\Oo) \subset {\cal N} \subset {\cal R}_{\omega}
(\Oo_1) \,. $$  
\\[6pt]
(e) Inner and outer regularity
\begin{equation}
 {\cal R}_{\omega}(\Oo) = \left( \bigcup_{\overline{\Oo_I} \subset \Oo}
 {\cal R}_{\omega}(\Oo_I) \right) '' =
 \bigcap_{\Oo_1 \supset \overline{\Oo}} {\cal R}_{\omega}(\Oo_1)
\end{equation}
holds for all regular diamonds $\Oo$.
\\[6pt]
(f) If $\omega$ is pure (an Hadamard vacuum), then we have Haag-Duality
$$ {\cal R}_{\omega}(\Oo)' = {\cal R}_{\omega}(\Oo^{\perp}) $$
for all regular diamonds $\Oo$. (By the same arguments as in {\rm
[65 (Prop.\ 6)]}, Haag-Duality extends to all pure (but not necessarily
quasifree or Hadamard) states $\omega$ which are locally normal
(hence, by (d), locally quasiequivalent) to any Hadamard vacuum.)
\\[6pt]
(g) Local primarity holds for all regular diamonds, that is, for
each regular diamond $\Oo$, ${\cal R}_{\omega}(\Oo)$ is a factor.
Moreover, ${\cal R}_{\omega}(\Oo)$ is isomorphic to the unique
hyperfinite type ${\rm III}_1$ factor if $\Oo^{\perp}$
is non-void. In this case, ${\cal R}_{\omega}(\Oo^{\perp})$ is
also hyperfinite and of type ${\rm III}_1$, and if $\omega$ is
pure, ${\cal R}_{\omega}(\Oo^{\perp})$ is again a factor.
Otherwise, if $\Oo^{\perp} = \emptyset$, then
${\cal R}_{\omega}(\Oo)$ is a type ${\rm I}_{\infty}$ factor.
\end{Theorem}
{\it Proof.} The key point in the proof is that, by results which
for the cases relevant here are to large extend due to Araki [1],
the above statement can be equivalently translated into statements
about the structure of the one-particle space, i.e.\ essentially the
symplectic space $(K,\kappa)$ equipped with the scalar product
$\lambda_{\omega}$. We shall use, however, the formalism of [40,45].
Following that, given a symplectic space $(K,\kappa)$ and
$\mu \in {\sf q}(K,\kappa)$ one calls a real linear map
${\bf k}: K \to H$ a {\it one-particle Hilbertspace structure} for
$\mu$ if (1) $H$ is a complex Hilbertspace, (2) the complex linear
span of ${\bf k}(K)$ is dense in $H$ and (3)
$$ \langle {\bf k}(x),{\bf k}(y) \rangle = \lambda_{\mu}(x,y)
= \mu(x,y) + \frac{i}{2}\kappa(x,y) $$
for all $x,y \in K$. It can then be shown (cf.\ [45 (Appendix A)])
that the GNS-representation of the quasifree state $\omega_{\mu}$
on $\AA[K,\kappa]$ may be realized in the following way:
$\Hh_{\omega_{\mu}} = F_s(H)$, the Bosonic Fock-space over the one-particle
space $H$, $\Omega_{\omega_{\mu}}$ = the Fock-vacuum, and
$$ \pi_{\omega_{\mu}}(W(x)) = {\rm e}^{i(a({\bf k}(x)) +
 a^+({\bf k}(x)))^-}\,, \quad
x \in K\, ,$$
where $a(\,.\,)$ and $a^+(\,.\,)$ are the Bosonic annihilation and
creation operators, respectively.

Now it is useful to define the symplectic complement
$F^{\tt v} := \{\chi \in H : {\sf Im}\,\langle \chi,\phi \rangle = 0
\ \ \forall \phi \in F \}$ for $F \subset H$, since it is known
that 
\begin{itemize}
\item[(i)] ${\cal R}_{\omega_{\mu}}(\Oo)$ is a factor \ \ \  iff\ \ \ 
  $ {\bf k}(K(\Oo))^- \cap {\bf k}(K(\Oo))^{\tt v} = \{0\}$,
\item[(ii)] ${\cal R}_{\omega_{\mu}}(\Oo)' = {\cal R}_{\omega_{\mu}}
(\Oo^{\perp})$\ \ \  iff\ \ \ 
 ${\bf k}(K(\Oo))^{\tt v} = {\bf k}(K(\Oo^{\perp}))^-$,
\item[(iii)] $\bigcap_{\Oo \supset C} {\cal R}_{\omega_{\mu}}(\Oo)
= \CC \cdot 1$ \ \ \  iff\ \ \ 
 $\bigcap_{\Oo \supset C}{\bf k}(K(\Oo))^- = \{0\}\,,$
\end{itemize}
cf.\ [1,21,35,49,58].

After these preparations we can commence with the proof of the various
statements of our Theorem.
\\[6pt]
(a) Let ${\bf k}: K \to H$ be the one-particle Hilbertspace structure
of $\omega$. The local one-particle spaces ${\bf k}(K(\Oo_G))^-$ of 
regular diamonds $\Oo_G$ based on $G \subset \Sigma$ are topologically
isomorphic to the completions of $C_0^{\infty}(G,\RR) \oplus
C_0^{\infty}(G,\RR)$ in the $H_{1/2} \oplus H_{-1/2}$-topology and
these are separable. Hence ${\bf k}(K)^-$, which is generated by a
countable set ${\bf k}(K(\Oo_{G_n}))$, for $G_n$ a sequence of
locally compact subsets of $\Sigma$ eventually exhausting $\Sigma$,
is also separable. The same holds then for the one-particle
Hilbertspace $H$ in which the complex span of ${\bf k}(K)$ is
dense, and thus separability is implied for $\Hh_{\omega} = F_s(H)$.
The infinite-dimensionality is clear.
\\[6pt]
(b) The local quasiequivalence has been proved in [66] and we refer to
that reference for further details. We just indicate that the 
proof makes use of the fact that the difference $\Lambda - \Lambda_1$
of the spatio-temporal two-point functions of any pair of
quasifree Hadamard states is on each causal normal neighbourhood
of any Cauchy-surface given by a smooth integral kernel ---
as can be directly read off from the Hadamard form --- and this turns
out to be sufficient for local quasiequivalence. The statement
about the unitary equivalence can be inferred from (g) below,
since it is known that every $*$-preserving isomorphism between
von Neumann algebras of type III acting on separable Hilbertspaces
is given by the adjoint action of a unitary operator which maps
the Hilbertspaces onto each other. See e.g.\ Thm.\ 7.2.9 and
Prop.\ 9.1.6 in [39]. 
\\[6pt]
(c) Here one uses that there exist Hadamard vacua, i.e.\ pure 
quasifree Hadamard states $\omega_{\mu}$. Since by Prop.\ 3.4
the topology of $\mu_{\Sigma}$ in $\DS$ is locally that of
$H_{1/2} \oplus H_{-1/2}$, one can show as in [66 (Chp.\ 4 and
Appendix)] that under the stated hypotheses about $C$ it holds
that $\bigcap_{\Oo \supset C} {\bf k}(K(\Oo))^- = \{0\}$ for the
one-particle Hilbertspace structures of Hadamard vacua. From
the local equivalence of the topologies induced by the dominating
scalar products of all quasifree Hadamard states (Prop.\ 3.4(e)),
this extends to the one-particle structures of all quasifree
Hadamard states. By (iii), this yields the statement (c).
\\[6pt]
(d) This is proved in [65] under the additional assumption that
the potential term $r$ is a positive constant. (The result was
formulated in [65] under the hypothesis that $\Sigma = \Sigma_1$,
but it is clear that the present statement without this hypothesis   
is an immediate generalization.) To obtain the general case
one needs in the spacetime deformation argument of [65]
the modification that the potential term $\hat{r}$ of the KG-field
on the new spacetime $(\hat{M},\hat{g})$ is equal to a positive
constant on its ultrastatic part while being equal to $r$ in a 
neighbourhood of $\Sigma$. We have used that procedure already in
the proof of Prop.\ 3.5, see also the proof of (f) below where 
precisely the said modification will be carried out in more detail.
\\[6pt]
(e) Inner regularity follows simply from the definition of the 
$\AA(\Oo)$; one deduces that for each $A \in \AA(\Oo)$ and each
$\epsilon > 0$ there exists some $\overline{\Oo_I} \subset \Oo$
and $A_{\epsilon} \in \AA(\Oo_I)$ so that
$||\,A - A_{\epsilon}\,|| < \epsilon$. It is easy to see that
inner regularity is a consequence of this property.

 So we focus now on the outer regularity.
 Let $\Oo = \Oo_G$ be based on the subset $G$ of the Cauchy-surface
$\Sigma$. Consider the symplectic space $(\DS,\delta_{\Sigma})$
and the dominating scalar product $\mu_{\Sigma}$ induced by $\mu
\in {\sf q}(\DS,\delta_{\Sigma})$, where $\omega_{\mu} = \omega$;
the corresponding one-particle Hilbertspace structure we denote by
${\bf k}_{\Sigma}: \DS \to H_{\Sigma}$. Then we denote by
${\cal W}({\bf k}_{\Sigma}({\cal D}_G))$ the von Neumann algebra in 
$B(F_s(H_{\Sigma}))$ generated by the unitary groups of the
operators $(a({\bf k}_{\Sigma}(\uu)) + a^+({\bf k}_{\Sigma}(\uu)))^-$
where $\uu$ ranges over ${\cal D}_G := C_0^{\infty}(G,\RR) \oplus
C_0^{\infty}(G,\RR)$. So ${\cal W}({\bf k}_{\Sigma}({\cal D}_G)) =
{\cal R}_{\omega}(\Oo_G)$. It holds
generally that $\bigcap_{G_1 \supset \overline{G}} {\cal W}({\bf k}_{\Sigma}
({\cal D}_{G_1})) = {\cal W}(\bigcap_{G_1 \supset \overline{G}}
{\bf k}_{\Sigma}({\cal D}_{G_1})^-)$ [1], hence, to establish outer 
regularity, we must show that
\begin{equation}
  \bigcap_{G_1 \supset \overline{G}} {\bf k}_{\Sigma}({\cal D}_{G_1})^-
 = {\bf k}_{\Sigma}({\cal D}_G)^-\,.
\end{equation}
In [65] we have proved that the ultrastatic vacuum $\omega^{\circ}$
of the KG-field with potential term $\equiv 1$ over the ultrastatic
spacetime $(M^{\circ},g^{\circ}) = (\RR \times \Sigma,dt^2 \oplus
(-\gamma))$ (where $\gamma$ is any complete Riemannian metric on
$\Sigma$) satisfies Haag-duality. That means, we have
\begin{equation}
{\cal R}^{\circ}_{\omega^{\circ}}(\Oo_{\circ})' =
{\cal R}^{\circ}_{\omega^{\circ}}(\Oo_{\circ}^{\perp})
\end{equation}
for any regular diamond $\Oo_{\circ}$ in $(M^{\circ},g^{\circ})$
which is based on any of the Cauchy-surfaces $\{t\}\times \Sigma$ in
the natural foliation, and we have put a ``$\circ$'' on the local
von Neumann algebras to indicate that they refer to a KG-field
over $(M^{\circ},g^{\circ})$. But since we have inner regularity
for ${\cal R}^{\circ}_{\omega^{\circ}}(\Oo_{\circ}^{\perp})$ ---
by the very definition --- the outer regularity of ${\cal R}^{\circ}
_{\omega^{\circ}}(\Oo_{\circ})$ follows from the Haag-duality (3.21).
Translated into conditions on the one-particle Hilbertspace
structure ${\bf k}^{\circ}_{\Sigma} : \DS \to H^{\circ}_{\Sigma}$
of $\omega^{\circ}$, this means that the equality
\begin{equation}
 \bigcap_{G_1 \supset \overline{G}} {\bf k}^{\circ}_{\Sigma}
({\cal D}_{G_1})^- = {\bf k}^{\circ}_{\Sigma}({\cal D}_G)^-
\end{equation}
holds. Now we know from Prop.\ 3.5 that $\mu_{\Sigma}$ induces
locally the $H_{1/2} \oplus H_{-1/2}$-topology on $\DS$. However,
this coincides with the topology locally induced by $\mu^{\circ}_{\Sigma}$
on $\DS$ (cf.\ (3.11)) --- even though $\mu^{\circ}_{\Sigma}$ may,
in general, not be viewed as corresponding to an Hadamard vacuum
of the KG-field over $(M,g)$. Thus the required relation (3.20)
is implied by (3.22).
\\[6pt]
(f) In view of outer regularity it is enough to show that, given
any $\Oo_1 \supset \overline{\Oo}$, it holds that 
\begin{equation}
 {\cal R}_{\omega}(\Oo^{\perp})' \subset {\cal R}_{\omega}(\Oo_1)\,.
\end{equation}
The demonstration of this property relies on a spacetime deformation
argument similar to that used in the proof of Prop.\ 3.5. Let
$G$ be the base of $\Oo$ on the Cauchy-surface $\Sigma$ in $(M,g)$.
Then, given any other open, relatively compact subset $G_1$ of
$\Sigma$ with  $\overline{G} \subset G_1$, we have shown in
[65] that there exists an ultrastatic spacetime $(\hat{M},\hat{g})$
with the properties (1) and (2) in the proof of Prop.\ 3.5, and with
the additional property that there is some $t < t_0$ such that
$$ \left( {\rm int}\,\hat{J}(G) \cap \Sigma_t \right )^- \subset
 {\rm int}\, \hat{D}(G_1) \cap \Sigma_t\,.$$
Here, $\Sigma_t = \{t\} \times \Sigma$ are the Cauchy-surfaces in the
natural foliation of the ultrastatic part of $(\hat{M},\hat{g})$.
The hats indicate that the causal set and the domain of dependence
are to be taken in $(\hat{M},\hat{g})$. This implies that we can find
some regular diamond $\Oo^t := {\rm int}\hat{D}(S^t)$ in
$(\hat{M},\hat{g})$ based on a subset $S^t$ of $\Sigma_t$ which
satisfies
\begin{equation}
\left( {\rm int}\, \hat{J}(G) \cap \Sigma_t \right)^-
\subset S^t \subset
{\rm int}\,\hat{D}(G_1) \cap \Sigma_t \,.
\end{equation}
Setting $\hat{\Oo} := {\rm int}\, \hat{D}(G)$ and
$\hat{\Oo}_1 := {\rm int}\,\hat{D}(G_1)$, one derives from (3.24)
the relations
\begin{equation}
 \hat{\Oo} \subset \Oo^t \subset \hat{\Oo}_1 
\,.
\end{equation}
These are equivalent to
\begin{equation}
\hat{\Oo}_1^{\perp} \subset (\Oo^t)^{\perp} \subset \hat{\Oo}^{\perp}
\end{equation}
where $\perp$ is the causal complementation in $(\hat{M},\hat{g})$.

Now as in the proof of Prop.\ 3.5, the given Hadamard vacuum $\omega$
on the Weyl-algebra $\AA[K,\kappa]$ of the KG-field over $(M,g)$
induces an Hadamard vacuum $\hat{\omega}$ on the Weyl-algebra 
$\AA[\hat{K},\hat{\kappa}]$ of the KG-field over $(\hat{M},\hat{g})$
whose potential term $\hat{r}$ is $1$ on the ultrastatic
part of $(\hat{M},\hat{g})$. Then by Prop.\ 6 in [65] we have
Haag-duality 
\begin{equation}
 \hat{\cal R}_{\hat{\omega}}(\hat{\Oo_t}^{\perp}) ' =
\hat{\cal R}_{\hat{\omega}}(\hat{\Oo_t})
\end{equation}
for all regular diamonds $\hat{\Oo_t}$ with base on
$\Sigma_t$; we have put hats on the von Neumann algebras
to indicate that they refer to $\AA[\hat{K},\hat{\kappa}]$.
(This was proved in [65] assuming that $(\hat{M},\hat{g})$
is globally ultrastatic. However, with the same argument, based on 
primitive causality, as we use it next to pass from (3.28) to  
(3.30), one can easily establish that (3.27) holds if only 
$\Sigma_t$ is, as here, a member in the natural foliation of the
ultrastatic part of $(\hat{M},\hat{g})$.)
Since $\Oo^t$ is a regular diamond based on $\Sigma_t$, we obtain
$$\hat{\cal R}_{\hat{\omega}}((\Oo^t)^{\perp})' =
\hat{\cal R}_{\hat{\omega}}(\Oo^t) $$
and thus, in view of (3.25) and (3.26),
\begin{equation}
\hat{\cal R}_{\hat{\omega}}(\hat{\Oo}^{\perp})'
\subset \hat{\cal R}_{\hat{\omega}}((\Oo^t)^{\perp})'
= \hat{\cal R}_{\hat{\omega}}(\Oo^t) \subset 
\hat{\cal R}_{\hat{\omega}}(\hat{\Oo}_1)\,.
\end{equation}
Now recall (see proof of Prop.\ 3.5) that $(\hat{M},\hat{g})$
coincides with $(M,g)$ on a causal normal neighbourhood $N$ of 
$\Sigma$. Primitive causality (Prop.\ 3.2) then entails 
\begin{equation}
\hat{\cal R}_{\hat{\omega}}(\hat{\Oo}^{\perp} \cap N)'
\subset \hat{\cal R}_{\hat{\omega}}(\hat{\Oo}_1 \cap N) \,.
\end{equation}
On the other hand, $\hat{\Oo}^{\perp} = {\rm int} \hat{D}(\Sigma
\backslash G)$ and $\hat{\Oo}_1$ are diamonds in $(\hat{M},\hat{g})$
based on $\Sigma$. Since $(M,g)$ and $(\hat{M},\hat{g})$
coincide on the causal normal neighbourhood $N$ of $\Sigma$,
one obtains that ${\rm int}\,D(\tilde{G}) \cap N =
{\rm int}\, \hat{D}(\tilde{G}) \cap N$ for all $\tilde{G} \in
\Sigma$. Hence, with $\Oo = {\rm int}\,D(G)$,
$\Oo_1 = {\rm int}\, D(G_1)$ (in $(M,g)$), we have that (3.23) entails
$$ {\cal R}_{\omega}(\Oo^{\perp} \cap N)' \subset
{\cal R}_{\omega}(\Oo_1 \cap N)
$$
(cf.\ the proof of Prop.\ 3.5) where the causal complement $\perp$
is now taken in $(M,g)$. Using primitive causality once more,
we deduce that 
\begin{equation}
 {\cal R}_{\omega}(\Oo^{\perp})' \subset {\cal R}_{\omega}(\Oo_1)\,.
\end{equation}
The open, relatively compact subset $G_1$ of $\Sigma$ was 
arbitrary up to the constraint $\overline{G} \subset G_1$.
Therefore, we arrive at the conclusion that the required inclusion
(3.23) holds of all $\Oo_1 \supset \overline{\Oo}$.
\\[6pt]
(g) Let $\Sigma$ be the Cauchy-surface on which $\Oo$ is based.
 For the local primarity one uses, as in (c), the existence of
Hadamard vacua $\omega_{\mu}$ and the fact (Prop.\ 3.5) that
$\mu_{\Sigma}$ induces locally the $H_{1/2} \oplus H_{-1/2}$-topology;
then one may use the arguments of [66 (Chp.\ 4 and Appendix)]
to show that due to the regularity of the boundary
$\partial G$ of the base $G$ of $\Oo$
there holds
$$ {\bf k}(K(\Oo))^- \cap {\bf k}(K(\Oo))^{\tt v} = \{ 0 \}$$
for the one-particle Hilbertspace structures of Hadamard vacua.
As in the proof of (c), this can be carried over to the
one-particle structures of all quasifree Hadamard states since
they induce locally on the one-particle spaces the same topology,
see [66 (Chp.\ 4)]. We note that for Hadamard vacua the local
primarity can also be established using (3.18) together with Haag-duality
and primitive causality purely at the algebraic level, without
having to appeal to the one-particle structures. 

The type ${\rm III}_1$-property of ${\cal R}_{\omega}(\Oo)$ is then
derived using Thm.\ 16.2.18 in [3] (see also [73]).
We note that for some points $p$ in the boundary $\partial G$ of $G$, $\Oo$
admits domains which are what is in Sect.\ 16.2.4 of [3] called
``$\beta_p$-causal sets'', as a consequence of the regularity of
$\partial G$ and the assumption $\Oo^{\perp} \neq \emptyset$.
We further note that it is straightforward to prove that
the quasifree Hadamard states of the KG-field over $(M,g)$
possess at each point in $M$ scaling limits (in the sense of
Sect.\ 16.2.4 in [3], see also [22,32]) which are equal to the
theory of the massless KG-field in Minkowski-spacetime. Together
with (a) and (c) of the present Theorem this shows that the 
the assumptions of Thm.\ 16.2.18 in [3] are fulfilled, and the
${\cal R}_{\omega}(\Oo)$ are type ${\rm III}_1$-factors for all
regular diamonds $\Oo$ with $\Oo^{\perp} \neq \emptyset$.
The hyperfiniteness follows from the split-property (d) and the
regularity (e), cf.\ Prop.\ 17.2.1 in [3]. The same arguments
may be applied to ${\cal R}_{\omega}(\Oo^{\perp})$, yielding
its type ${\rm III}_1$-property (meaning that in its central
decomposition only type ${\rm III}_1$-factors occur) and
hyperfiniteness. If $\omega$ is an Hadamard vacuum, then
${\cal R}_{\omega}(\Oo^{\perp}) = {\cal R}_{\omega}(\Oo)'$ is
a factor unitarily equivalent to ${\cal R}_{\omega}(\Oo)$.
For the last statement note that $\Oo^{\perp} = \emptyset$ implies
that the spacetime has a compact Cauchy-surface on which $\Oo$
is based. In this case ${\cal R}_{\omega}(\Oo) =
\pi_{\omega}(\AA[K,\kappa])''$ (use the regularity of $\partial G$,
and (c), (e) and primitive causality). But since $\omega$ is
quasiequivalent to any Hadamard vacuum by the relative compactness of 
$\Oo$, ${\cal R}_{\omega}(\Oo) = \pi_{\omega}(\AA[K,\kappa])''$
is a type ${\rm I}_{\infty}$-factor. $\Box$
\\[10pt]
We end this section and therefore, this work, with a few concluding
remarks.
 
First we note that the split-property signifies a strong notion of
statistical independence. It can be deduced from constraints on the
phase-space behaviour (``nuclearity'') of the considered quantum
field theory. We refer to [9,31] for further information and
also to [62] for a review, as a discussion of these issues lies 
beyond the scope of of this article. The same applies to a discussion
of the property of the local von Neumann algebras ${\cal R}_{\omega}(\Oo)$
to be hyperfinite and of type ${\rm III}_1$. We only mention that
for quantum field theories on Minkowski spacetime it can be established
under very general (model-independent)
 conditions that the local (von Neumann) observable
algebras are hyperfinite and of type ${\rm III}_1$, and refer the reader to
[7] and references cited therein. However, the property of the
local von Neumann algebras to be of type ${\rm III}_1$, together with
the separability of the GNS-Hilbertspace $\Hh_{\omega}$, has an
important consequence which we would like to point out (we have
used it implicitly already in the proof of Thm.\ 3.6(b)):
$\Hh_{\omega}$ contains a dense subset ${\sf ts}(\Hh_{\omega})$ of vectors
which are cyclic and separating for all ${\cal R}_{\omega}(\Oo)$
whenever $\Oo$ is a diamond with $\Oo^{\perp} \neq \emptyset$. But
so far it has only been established in special cases that $\Omega_{\omega}
\in {\sf ts}(\Hh_{\omega})$, see [64]. At any rate, when
$\Omega \in {\sf ts}(\Hh_{\omega})$ one may consider for a pair of
regular diamonds $\Oo_1,\Oo_2$ with $\overline{\Oo_1} \subset \Oo_2$
and $\Oo_2^{\perp}$ nonvoid the modular operator $\Delta_2$
of ${\cal R}_{\omega}(\Oo_2)$,$\Omega$ (cf.\ [39]). The split property
and the factoriality of ${\cal R}_{\omega}(\Oo_1)$ and ${\cal R}_{\omega}
(\Oo_2)$ imply the that the map
\begin{equation}
 \Xi_{1,2} : A \mapsto \Delta^{1/4}_2 A \Omega\,, \quad A \in 
{\cal R}_{\omega}(\Oo_1)\,,
\end{equation}
is compact [8]. As explained in [8],
``modular compactness'' or ``modular nuclearity'' may be viewed
as suitable generalizations of ``energy compactness'' or
``energy nuclearity'' to curved spacetimes as notions to measure
the phase-space behaviour of a quantum field theory
(see also [65]). Thus an interesting
question would be if the maps (3.31) are even nuclear.

Summarizing it can be said that Thm.\ 3.6 shows that the nets of von
Neumann observable algebras of the KG-field over a globally hyperbolic
spacetime in the representations of quasifree Hadamard states have
all the properties one would expect for physically reasonable
representations. This supports the point of view that quasifree
Hadamard states appear to be a good choice for physical states
of the KG-field over a globally hyperbolic spacetime. Similar results
are expected to hold also for other linear fields. 

Finally, the reader will have noticed that we have been considering
exclusively the quantum theory of a KG-field on a {\it globally hyperbolic}
spacetime. For recent developments concerning quantum fields in the
background of non-globally hyperbolic spacetimes, we
refer to [44] and references cited there.
\\[24pt]
{\bf Acknowledgements.} I would like to thank D.\ Buchholz for 
valueable comments on a very early draft of Chapter 2. Moreover,
I would like to thank C.\ D'Antoni, R.\ Longo,  J.\ Roberts
and L.\ Zsido for
their hospitiality, and their interest in quantum field theory
in curved spacetimes. I also appreciated conversations with R.\ Conti,
D.\ Guido and L.\ Tuset on various parts of the material of the
present work. 
\\[28pt]
\noindent
{\Large {\bf Appendix}}
\\[24pt]
{\bf Appendix A} 
\\[18pt]
For the sake of completeness, we include here the interpolation argument
in the form we use it in the proof of Theorem 2.2 and in Appendix B
below. It is a standard argument based on Hadamard's three-line-theorem,
cf.\ Chapter IX in [57].
\\[10pt]
{\bf Lemma A.1}
{\it 
Let $\Ff,\Hh$ be  complex Hilbertspaces, $X$ and $Y$ two non-negative,
injective, selfadjoint operators  in $\Ff$ and $\Hh$, respectively,
 and $Q$ a bounded linear operator  $\Hh \to \Ff$
such that $Q{\rm Ran}(Y) \subset {\rm dom}(X)$. 
 Suppose that the operator $XQY$ admits a
bounded extension $T :\Hh \to \Ff$. Then for all $0 \leq \tau \leq 1$,
it holds that $Q{\rm Ran}(Y^{\tau}) \subset {\rm dom}(X^{\tau})$, and
the operators $X^{\tau}QY^{\tau}$ are bounded by $||\,T\,||^{\tau}
||\,Q\,||^{1 - \tau}$.   }
\\[10pt]
{\it Proof.} The operators $\ln(X)$ and $\ln(Y)$ are (densely defined)
selfadjoint operators. Let the vectors $x$ and $y$ belong to the
spectral subspaces of $\ln(X)$ and $\ln(Y)$, respectively, corresponding
to an arbitrary finite intervall. Then the functions
$\CC \owns z \mapsto {\rm e}^{z\ln(X)}x$ and
$\CC \owns z \mapsto {\rm e}^{z\ln(Y)}y$ are holomorphic. Moreover,
${\rm e}^{\tau \ln(X)}x = X^{\tau}x$ and ${\rm e}^{\tau \ln(Y)}y =
Y^{\tau}y$ for all real $\tau$. Consider the function
$$ F(z) := \langle {\rm e}^{\overline{z}\ln(X)}x,Q{\rm e}^{z\ln(Y)}y
\rangle_{\Ff} \,.$$
It is easy to see that this function is holomorphic on $\CC$, and
also that the function is uniformly bounded for $z$ in the
strip $\{z : 0 \leq {\sf Re}\,z \leq 1 \}$.
For $z = 1 + it$, $t \in \RR$, one has
$$ |F(z)| = |\langle {\rm e}^{-it\ln(X)}x,XQY{\rm e}^{it\ln(Y)}y \rangle_{\Ff} |
  \leq ||\,T\,||\,||\,x\,||_{\Ff}||\,y\,||_{\Hh} \,,$$
and for $z = it$, $t \in \RR$,
$$ |F(z)| = |\langle {\rm e}^{-it\ln(X)}x,Q{\rm e}^{it\ln(Y)}y \rangle_{\Ff} |
  \leq ||\,Q\,||\,||\,x\,||_{\Ff}||\,y\,||_{\Hh} \,.$$
By Hadamard's three-line-theorem, it follows that for all $z = \tau + it$
in the said strip there holds the bound
$$ |F(\tau + it)| \leq ||\,T\,||^{\tau}||\,Q\,||^{1 - \tau}||\,x\,||_{\Ff}
||\,y\,||_{\Hh}\,.$$
As $x$ and $y$ were arbitrary members of the finite spectral intervall
subspaces, the last estimate extends to all $x$ and $y$ lying in cores
for the operators
  $X^{\tau}$ and $Y^{\tau}$, from which the 
the claimed statement follows. $\Box$
\\[24pt]
{\bf Appendix B}
\\[18pt]
For the convenience of the reader we collect here two well-known
results about Sobolev norms on manifolds which are used in
the proof of Proposition 3.5. The notation is as follows.
$\Sigma$ and $\Sigma'$ will denote smooth, finite dimensional manifolds
(connected, paracompact, Hausdorff); $\gamma$ and $\gamma'$ are
complete Riemannian metrics on $\Sigma$ and $\Sigma'$, respectively.
Their induced volume measures are denoted by $d\eta$ and $d\eta'$.
We abbreviate by $A_{\gamma}$ the selfadjoint extension
in $L^2(\Sigma,d\eta)$
of the operator $-\Delta_{\gamma} +1$ on $C_0^{\infty}(\Sigma)$,
where $\Delta_{\gamma}$ is the Laplace-Beltrami operator on $(\Sigma,\gamma)$;
note that [10] contains a proof that $(-\Delta_{\gamma} + 1)^k$
is essentially selfadjoint on $C_0^{\infty}(\Sigma)$ for all $k \in \NN$.
$A'$ will be defined similarly with respect to the corresponding
objects of $(\Sigma',\gamma')$. As in the main text, the
$m$-th Sobolev scalar product is $\langle u,v \rangle_{\gamma,m}
= \langle u,A_{\gamma}^{m}v \rangle$ for $u,v \in C_0^{\infty}(\Sigma)$ and
$m \in \RR$, where $\langle\,.\,,\,.\, \rangle$ is the scalar product
of $L^2(\Sigma,d\eta)$. Anagolously we define $\langle\,.\,,\,.\,\rangle_{
\gamma',m}$. For the corresponding norms we write $||\,.\,||_{\gamma,m}$,
resp., $||\,.\,||_{\gamma',m}$.
\\[10pt]
{\bf Lemma B.1} 
{\it    (a) Let $\chi \in C_0^{\infty}(\Sigma)$. Then there is for each
$m \in \RR$ a constant $c_m$ so that
$$ ||\,\chi u \,||_{\gamma,m} \leq c_m  ||\, u \,||_{\gamma,m}\,,
\quad u \in C_0^{\infty}(\Sigma) \,.$$ 
\\[6pt]
(b) Let $\phi \in C^{\infty}(\Sigma)$ be strictly positive and 
$G \subset \Sigma$ open and relatively compact. Then there are 
for each $m \in \RR$ two positive constants $\beta_1,\beta_2$ so that
$$ \beta_1||\,\phi u\,||_{\gamma,m} \leq ||\,u\,||_{\gamma,m}
\leq \beta_2||\,\phi u\,||_{\gamma,m}\,, \quad u \in C_0^{\infty}(G)\,.$$ }
{\it Proof.} (a) We may suppose that $\chi$ is real-valued (otherwise
we treat real and imaginary parts separately). A tedious but straightforward
calculation shows that the claimed estimate is fulfilled for all $m =2k$,
$k \in \NN_0$. Hence $A^k \chi A^{-k}$ extends to a bounded operator
on $L^2(\Sigma,d\eta)$, and the same is true of the adjoint
$A^{-k}\chi A^k$. 
 Thus by the interpolation argument, cf.\ Lemma A.1,
$A^{\tau k} \chi A^{-\tau k}$ is bounded for all $-1 \leq \tau \leq 1$.
This yields the stated estimate.
\\[6pt]
(b) This is a simple corollary of (a). For the first estimate, note
that we may replace $\phi$ by a smooth function with compact support.
Then note that the second estimate is equivalent to 
$||\,\phi^{-1}v\,||_{\gamma,m} \leq \beta_2||\,v\,||_{\gamma,m}$,
$v \in C_0^{\infty}(G)$, and again we use that instead of $\phi^{-1}$
we may take a smooth function of compact support. $\Box$
\\[10pt]
{\bf Lemma B.2} 
{\it
     Let $(\Sigma,\gamma)$ and $(\Sigma',\gamma')$ be two complete 
Riemannian manifolds, $N$ and $N'$ two open subsets of $\Sigma$ and
$\Sigma'$, respectively, and $\Psi : N \to N'$ a diffeomorphism.
Given $m \in \RR$ and some open, relatively compact subset $G$ of
$\Sigma$ with $\overline{G} \subset N$, there are two positive
constants $b_1,b_2$ such that
$$ b_1||\,u\,||_{\gamma,m} \leq ||\,\Psi^*u\,||_{\gamma',m}
\leq b_2||\,u\,||_{\gamma,m} \,, \quad u \in C_0^{\infty}(G)\,,$$
where $\Psi^*u := u \lcrc \Psi^{-1}$. }
\\[10pt]
{\it Proof.} Again it is elementary to check that such a result
is true for $m = 2k$ with $k \in \NN_0$. One infers that, choosing
$\chi \in C_0^{\infty}(N)$ with $\chi|G \equiv 1$ and setting
$\chi' := \Psi^*\chi$, there is for each $k \in \NN_0$ a
positive constant $b$ fulfilling 
$$ ||\,A^k\chi\Psi_*\chi'v\,||_{\gamma,0} \leq
  b\,||\,(A')^kv\,||_{\gamma',0}\,, \quad v \in C_0^{\infty}(\Sigma')\,;$$
here $\Psi_*v := v \lcrc \Psi$. Therefore,
$$ A^k\lcrc \chi\lcrc \Psi_*\lcrc\chi'\lcrc(A')^{-k} $$
extends to a bounded operator $L^2(\Sigma',d\eta') \to L^2(\Sigma,d\eta)$
for each $k \in \NN_0$. Interchanging the roles of $A$ and $A'$, one
obtains that also 
$$ (A')^k\lcrc \chi'\lcrc \Psi^*\lcrc \chi\lcrc A^{-k} $$
extends, for each $k \in \NN_0$, to a bounded operator
$L^2(\Sigma,d\eta) \to L^2(\Sigma',d\eta')$. The boundedness transfers
to the adjoints of these two operators. Observe then that for
$(\Psi_*)^{\dagger}$, the adjoint of $\Psi_*$, we have
$(\Psi_*)^{\dagger} = \rho^2\lcrc(\Psi^*)$ on $C_0^{\infty}(N)$, and
similarly, for the adjoint $(\Psi^*)^{\dagger}$ of $\Psi^*$ we have
$(\Psi^*)^{\dagger} = \Psi_* \lcrc \rho^{-2}$ on $C_0^{\infty}(N')$,
where $\rho^2 = \Psi^*d\eta/d\eta'$ is a smooth density function
on $N'$, cf.\ eqn.\ (3.14).
It can now easily be worked out that the interpolation argument
of Lemma A.1 yields again the claimed result.
\begin{flushright}
 $\Box$
\end{flushright}
{\small  
 }
\end{document}